%% file: main.tex
\documentclass[10pt,reqno]{article}
\usepackage[a4paper,margin=1.1in]{geometry}

\usepackage[T1]{fontenc}
\usepackage[utf8]{inputenc}
\usepackage{lmodern,microtype}
\usepackage[parfill,indent]{parskip}
\usepackage{bbm}
\usepackage[cal=boondoxo]{mathalfa}
\usepackage{xspace}
\usepackage{setspace}
\usepackage{todonotes}

\usepackage{amsmath,amssymb,amsthm,mathtools}
\usepackage[authoryear,round]{natbib}
\usepackage{xcolor}
\usepackage[shortlabels]{enumitem}
\setlist[itemize]{noitemsep}
\usepackage{hyperref}

\hypersetup{
    colorlinks=true,
    breaklinks=true,
    linkcolor=blue,
    urlcolor=blue,
    anchorcolor=blue,
    citecolor=blue!50!black,
    pdftitle={Equilibrium measures for anisotropic exponential interaction
              energies with applications to transformers},
    pdfauthor={Oleksii Kachaiev, Martin Burger, Samira Kabri,
               Yury Korolev, Tim Roith},
    pdfsubject={Equilibrium measures for anisotropic exponential
                interaction energies and mean-field models of self-attention},
    pdfkeywords={equilibrium measures, interaction energies,
                 self-attention, transformers, Brownian motion,
                 Feynman--Kac},
    pdfdisplaydoctitle=true,
    bookmarksnumbered=true
}

\usepackage{bookmark}
\usepackage{algorithm}
\usepackage[noend]{algpseudocode}
\usepackage{subcaption}
\setlist[enumerate]{itemsep=0.2em,topsep=0.4em}

\newtheorem{theorem}{Theorem}

\newtheorem{proposition}[theorem]{Proposition}
\newtheorem{corollary}[theorem]{Corollary}

\newtheorem*{conjecture*}{Conjecture}
\theoremstyle{remark}
\newtheorem{remark}{Remark}
\newcommand{\bbR}{\mathbb{R}}
\newcommand{\bbN}{\mathbb{N}}
\newcommand{\bbS}{\mathbb{S}}
\newcommand{\bbE}{\mathbb{E}}
\newcommand{\bbP}{\mathbb{P}}
\newcommand{\cV}{\mathcal{V}}
\newcommand{\cE}{\mathcal{E}}
\newcommand{\cP}{\mathcal{P}}
\newcommand{\cH}{\mathcal{H}}
\newcommand{\cG}{\mathcal{G}}

\newcommand{\cL}{\mathcal{L}}
\newcommand{\cF}{\mathcal{F}}
\newcommand{\cA}{\mathcal{A}}
\newcommand{\cN}{\mathcal{N}}

\newcommand{\cZ}{\mathcal{Z}}
\newcommand{\Sph}{\bbS^{d-1}}
\newcommand{\Ball}{\mathbb B^d}
\newcommand{\norm}[1]{\left\|#1\right\|}
\newcommand{\abs}[1]{\left|#1\right|}
\newcommand{\inner}[2]{\left\langle#1,#2\right\rangle}

\newcommand{\dd}{\mathop{}\!\mathrm{d}}

\newcommand{\Tr}{\operatorname{Tr}}
\newcommand{\diag}{\operatorname{diag}}
\newcommand{\supp}{\operatorname{supp}}
\newcommand{\MMD}{\operatorname{MMD}}
\newcommand{\op}{\mathrm{op}}
\newcommand{\TV}{\mathrm{TV}}

\newcommand{\Id}{I}
\newcommand{\Exp}{\mathrm{Exp}}
\DeclareMathOperator{\diver}{div}

\newcommand{\weakstarto}{\rightharpoonup}
\newcommand{\one}{\mathbf{1}}

\definecolor{fuzzywuzzy}{rgb}{0.8, 0.4, 0.4}

\definecolor{trr}{rgb}{0.1, 0.7, 0.1}

\title{Equilibrium measures for anisotropic exponential interaction\\energies with applications to transformers}

\makeatletter
\renewcommand\@fnsymbol[1]{\@arabic{#1}}
\makeatletter
\author{%
  Oleksii Kachaiev%
  \thanks{MaLGa, Dipartimento di Matematica,
  Universit\`a degli Studi di Genova, Genoa, Italy.}\and%
  Martin Burger\thanks{Helmholtz Imaging, Deutsches Elektronen-Synchrotron DESY, Notkestr. 85,  22607 Hamburg, Germany.} \textsuperscript{,}\thanks{Department of Mathematics, University of Hamburg, Bundesstr. 55, 20146 Hamburg, Germany.}\and%
  Samira Kabri\footnotemark[2] \textsuperscript{,}\footnotemark[3] \and%
  Yury Korolev%
  \thanks{Department of Mathematical Sciences,
  University of Bath, Bath, UK.}\and%
  Tim Roith\thanks{School of Computation, Information and Technology,
  Technical University of Munich, Germany.} \textsuperscript{,}\thanks{Munich Center for Machine Learning, Munich, Germany.}\\
  {\small
    Correspondence to:
    \texttt{oleksii.kachaiev@edu.unige.it}%
  }  
}

\date{}

\begin{document}

\maketitle

\vspace{-15pt}

\begin{abstract}
We study the minimization of anisotropic exponential interaction
energies over probability measures on the sphere, motivated by
mean-field models of self-attention in transformer architectures. For every symmetric positive definite interaction matrix, we relate the minimizer to a corresponding anisotropic Brownian motion and use the Feynman--Kac formula and elliptic regularity theory to show that it admits a smooth, strictly positive density on the sphere. We further show that it is the only stationary measure with full support for the dynamics considered. We use the stochastic viewpoint to derive an algorithm for sampling from the minimizer. In dimension two, we also propose a numerical
scheme for approximating its density based on a Mathieu series expansion.
\end{abstract}

\section{Introduction}

The equilibrium behavior of interacting particle systems is a classical topic in statistical physics, including a variety of different phenomena such as clustering, screening, or crystallization (cf. \cite{tindemans1974exact,saff1997distributing,hardin2005minimal,serfaty2014ginzburg}). Analytical insights for large particle systems can be gained from a mean-field perspective, which opens the way to methods in calculus of variations (in the space of measures) or partial differential equations (cf. \cite{mackowiak1982mean,burger2008large,burger2014stationary,simione2015existence,carrillo2017geometry,canizo2024discrete,shu2024wasserstein}), recently also extending the range of applications far beyond physics, e.g. to systems in biology or social behavior. A recent unexpected development in the field is the study of large-scale behavior in transformer architectures, widely used in modern machine learning, through interacting particle systems. This perspective helps explain some of the rich behavior observed in computations. The normalization mechanisms in these architectures  naturally lead to problems on the sphere, whereas most of particle systems in other settings are studied in the full space. Moreover, even with a strong level of simplification relative to realistic transformer architectures, the anisotropy of the interaction remains important, distinguishing these problems the predominant isotropic settings in the study of equilibrium measures. 

In this paper, our main object of interest is the following energy functional on the space of Borel probability measures on the unit sphere \(\cP(\Sph)\), $d \geq 2$,
\begin{equation}
\label{eq:energy}
\cE_D(\mu):=\frac12
\iint_{\Sph\times\Sph}\! K_D(x,y) \,\mu(\dd x)\mu(\dd y),
\qquad \mu\in\cP(\Sph),
\end{equation}
where the interaction kernel is given by \( K_D(x,y) :=  e^{x^\top Dy} \)  with a symmetric and positive definite matrix \(D\in\bbR^{d\times d}\).

Variational problems of this form, with various interaction kernels $K_D$, arise
in several areas, including aggregation and swarming models, where attractive and repulsive interactions determine the shape of equilibrium configurations, e.g.~\citet{carrillo2014global}, 
models of collective alignment on the
sphere, e.g.~\citet{frouvelle2012dynamics}, and in particular mean-field self-attention models in transformer architectures, e.g.~\citet{sander2022sinkformers,geshkovski2025perspective} and subsequent works (due to the large and ever expanding literature on transformers we do not attempt a review; we restrict attention to works most directly relevant to our setting). For \(D=\lambda I\) with \( \lambda > 0 \), it is known~\citep{borodachov2019discrete} that the unique minimizer is the uniform distribution on the sphere; this has also been noted by~\citet{geshkovski2025perspective}  in the context of transformers.
In the anisotropic case, the properties of the minimizer depend heavily on the eigenvalues of the matrix $D$~\citep{burger2025analysis}. When $D$ is not positive definite, the minimizer may be supported on a proper subspace spanned by eigenvectors of $D$, and is therefore singular with respect to the surface measure. For $D\succ0$, however, it was conjectured that the minimizer has full support.
In this paper we prove this conjecture using techniques from stochastic analysis and elliptic PDE theory.

\paragraph{Main contributions and structure of the paper.} 
Our main contributions are as follows. \begin{itemize}
    \item We show that the minimizer of~\eqref{eq:energy} is the exit distribution of a certain stochastic process (killed anisotropic Brownian motion; Section~\ref{sec:probability}) and, furthermore, that it has a smooth and strictly positive density with respect to the surface measure on the sphere (Section~\ref{sec:density}), thus resolving the full-support conjecture of \citet{burger2025analysis}.
    \item We show that, for the dynamics considered, the minimizer is the only stationary measure having full support and that every other stationary measure is singular with respect to the surface measure (Section~\ref{sec:stationary-points-singular}).
    \item We analyze the stability of the minimizer (in the maximum mean discrepancy induced by the interaction kernel) with respect to perturbations in the matrix $D$ (Section~\ref{sec:stability}).
    \item We propose an algorithm for sampling from the minimizer using the stochastic process introduced in Section~\ref{sec:probability} (Section~\ref{subsec:comp-wos}).
    \item For $d=2$, we propose a numerical scheme for approximating the minimizing density based on a Mathieu series expansion (Section~\ref{sec:Mathieu}).
\end{itemize}

In Section~\ref{sec:var-prob} we collect the properties of the variational problem that will be used in the rest of the paper. Our notation is summarized below.

\paragraph{Notation.}
Throughout the paper, \(D\in\bbR^{d\times d}\) with \(d\geq2\) is a symmetric positive definite matrix, i.e. \( D^{\top} = D \) and \( D \succ 0 \).
For $y\in\bbR^d$ we denote by $\diag(y)\in\bbR^{d\times d}$ the diagonal matrix with $\diag(y)_{ii} = y_i$ for all $i\in\{1,\ldots, d\}$.
We use \(\norm{\cdot}\) for the Euclidean norm and \(\Sph=\{x\in\bbR^d:\norm{x}=1\}\) and \(\Ball:=\{y\in\bbR^d:\norm{y}<1\}\) for the unit sphere and the open unit ball, respectively. For matrices, we denote by \(I\), \(\det\), \(\Tr\), and \(\norm{\cdot}_{\op}\) respectively the identity, determinant, trace, and  operator norm induced by~\(\norm{\cdot}\).
We write \(a\wedge b:=\min\{a,b\}\). For a twice continuously differentiable function \(u\), we denote by \( \nabla^2 u \) its Hessian. Let \(\nabla_{\Sph}\) and \(\diver_{\Sph}\) denote the intrinsic
gradient and divergence.

We denote by \(\cP(\Sph)\) the space of Borel probability measures over \(\Sph\), by $\dd S$ the surface measure, and by \(\dd \sigma_{d-1}\)  the normalized surface measure.
We write \(\delta_z\) for the Dirac measure at \(z\),
\(\supp\mu\) for the support of \(\mu\), and \(F_\#\mu\) for its
push-forward under \(F\), defined by
\((F_\#\mu)(A)=\mu(F^{-1}(A))\).
We write \( \mu \perp \nu \) when the measures \( \mu \) and \( \nu \) are
mutually singular, that is, concentrated on disjoint Borel sets.
The spaces of real-valued continuous and smooth functions on \(\Sph\)
are denoted by \(C(\Sph)\) and \(C^\infty(\Sph)\), respectively.
We write \( C_c^2(\bbR^d) \) for the space of twice continuously differentiable functions on \( \bbR^d \) with compact support. Furthermore, we denote by $L^2(\Sph)$ the Hilbert space of Lebesgue functions on $\Sph$ with respect to the surface measure $\dd S$, equipped with the inner product
\begin{align*}
\langle f, g\rangle_{L^2}
:=  \int_{\Sph} f g \, \dd S,
\qquad f,g\in L^2(\Sph).
\end{align*}
For a sequence \((\mu_n)_{n\geq1}\subset\cP(\Sph)\) and
\(\mu\in\cP(\Sph)\), weak convergence, written \(\mu_n\weakstarto\mu\), means that
\begin{equation*}
\int_{\Sph} \varphi\,\dd\mu_n
\longrightarrow \int_{\Sph} \varphi\,\dd\mu
\qquad\text{for every } \varphi\in C(\Sph).
\end{equation*}
For a finite signed Borel measure \(\sigma\) on \(\Sph\), we denote its
total variation norm by
\[
\norm{\sigma}_{\TV}
:= \sup_{\substack{f\in C(\Sph)\\ \norm{f}_\infty\leq 1}}
\left|\int_{\Sph} f\,\dd\sigma\right|.
\]
The Wasserstein distance (of order $2$) is defined by
\[
W_2(\mu,\nu)
:=
\left(
\inf_{\gamma\in\Pi(\mu,\nu)}
\int_{\Sph\times\Sph}d_{\Sph}(x,y)^2\,\gamma(\dd x,\dd y)
\right)^{1/2},
\qquad 
\mu,\nu\in\cP(\Sph),
\]
where \(\Pi(\mu,\nu)\) is the set of Borel probability measures
on \(\Sph\times\Sph\) with marginals \(\mu\) and \(\nu\), and
\(d_{\Sph}(x,y):=\arccos(x^\top y)\) is the geodesic distance.
Since \(\Sph\) is compact, the Wasserstein distance metrizes the weak topology \citep[Corollary~6.13]{villani2009optimal}.

\section{Properties of the variational problem}
\label{sec:var-prob}

In this section we collect some results about the variational problem associated with the energy
\(\cE_D\) that will be used throughout the paper. We consider
\begin{equation}
\label{eq:variational-problem}
\underset{\mu\in\cP(\Sph)}{\operatorname{argmin}}\,
\left\{ \cE_D(\mu)=\frac12 \iint_{\Sph\times\Sph}K_{D}(x,y)\,\mu(\dd x)\mu(\dd y) \right\}
\qquad \text{with} \quad
K_{D}(x, y) = e^{x^\top Dy},
\end{equation}
where, as before, \(D\) is symmetric and positive definite.
We first establish existence and uniqueness of a global minimizer.
We then introduce the Euler--Lagrange
criterion, which will guide the construction of the minimizer in
Section~\ref{sec:equilibrium}. Finally, we relate the energy gap to a
distance between probability measures induced by the interaction kernel.

\begin{proposition}[Existence and uniqueness]
\label{prop:unique-minimizer}
The problem~\eqref{eq:variational-problem} admits a unique minimizer,
which we denote by \(\mu_D\).
\end{proposition}

\begin{proof}
The space \(\cP(\Sph)\) is weakly compact, and continuity of
\(K_D\) on \(\Sph\times\Sph\) makes \(\cE_D\) weakly
continuous. Hence the minimum is attained.
By \citet[Lemma~3.6]{burger2025analysis}, the energy \(\cE_D\) is
strictly convex on \(\cP(\Sph)\), that is 
\begin{equation}
\label{eq:convexity}
\cE_D((1-t)\mu+t\nu)
<(1-t)\cE_D(\mu)+t\cE_D(\nu),
\qquad \mu,\nu\in\cP(\Sph), \; \mu \neq\nu.
\end{equation}
Therefore the minimizer is unique.
\end{proof}

\paragraph{Euler--Lagrange optimality criterion.} For a finite signed Borel measure \(\rho\) on \(\Sph\), define its \emph{potential} 
by 
\begin{equation}
\label{eq:potential}
U_\rho(x)
:=\int_{\Sph}K_D(x,y)\,\rho(\dd y),
\qquad x\in\Sph,
\end{equation}
where we simplify notation by suppressing the dependence on $D$.

\begin{proposition}[Euler--Lagrange criterion  {\citep[Theorem~2.14 and Lemma~3.3]{bilyk2022positive}}]
\label{prop:euler-lagrange}
For \(\mu\in\cP(\Sph)\),
\begin{equation}
\label{eq:el}
\mu=\mu_D
\quad\Longleftrightarrow\quad
\begin{cases}
U_\mu(x)\geq2\cE_D(\mu)
&\text{for every }x\in\Sph,\\
U_\mu(x)=2\cE_D(\mu)
&\text{for every }x\in\supp\mu.
\end{cases}
\end{equation}
In particular, by uniqueness of the minimizer, any probability measure whose potential is constant
on \(\Sph\) is equal to \(\mu_D\).
\end{proposition}

The criterion says that the support of an equilibrium measure lies
where its potential is smallest. A measure with constant potential
on the whole sphere therefore provides a particularly convenient
candidate: once it is known to be a probability measure, minimality
follows. This observation will be the starting point of the
proofs of main results in Section~\ref{sec:equilibrium}.

\paragraph{Interpretation using reproducing kernel Hilbert spaces.}

We next show that the interaction kernel \(K_D\) is strictly positive
definite whenever \(D\) is symmetric and positive definite. Set
\[
w(x):=e^{x^\top Dx/2},
\qquad
g(u,v):=e^{-\norm{u-v}^2/2}.
\]
The identity
\begin{equation}
\label{eq:gaussian-factorization}
K_D(x,y)=w(x)w(y)g(D^{1/2}x,D^{1/2}y)
\end{equation}
shows that \(K_D\) is strictly positive definite:
the Gaussian kernel \(g\) is strictly positive definite
\citep[Proposition~5]{sriperumbudur2011universality},
the map \(x\mapsto D^{1/2}x\) is injective, and \(w\) is
strictly positive.

Thus there exists a unique
reproducing kernel Hilbert spaces (RKHS) \(\cH_D\) with reproducing
kernel \(K_D\)~\citep{aronszajn1950theory}, whose \emph{canonical feature map} is
\[
x\longmapsto K_D(x,\cdot)\in\cH_D,
\qquad x \in \Sph.
\]
For a finite signed Borel measure \(\rho\), its potential
\begin{equation}
\label{eq:u-rho-rkhs}
U_\rho
=\int_{\Sph}K_D(y,\cdot)\,\rho(\dd y)
\in\cH_D
\end{equation}
is its \emph{kernel mean embedding}~\citep[Definition~1]{simongabriel2018kernel}.
The Bochner integral~\eqref{eq:u-rho-rkhs} is well defined since
\(\norm{K_D(y,\cdot)}_{\cH_D}^2=K_D(y,y)\) is bounded on
\(\Sph\).
The reproducing property gives, for finite signed Borel measures
\(\rho,\eta\),
\begin{equation}
\label{eq:rkhs-inner-u}
\inner{U_\rho}{U_\eta}_{\cH_D}
=\iint_{\Sph\times\Sph}
K_D(x,y)\,\rho(\dd x)\eta(\dd y),
\end{equation}
so that
\begin{equation}
\label{eq:rkhs-norm-sq}
\cE_D(\mu) = \frac{1}{2} \norm{U_\mu}_{\cH_D}^2,
\qquad \mu \in \cP(\Sph).
\end{equation}

For \(\mu,\nu\in\cP(\Sph)\), the associated
\emph{maximum mean discrepancy}
\citep[Section~2.2]{gretton2012kernel} satisfies
\begin{align*}
\MMD_{D}(\mu,\nu)
&:=\sup_{\norm{f}_{\cH_D} \leq1}
\left|\int_{\Sph}f\,\dd(\mu-\nu)\right|
=\norm{U_\mu-U_\nu}_{\cH_D}.
\end{align*}

\begin{proposition}
\label{prop:kernel-norm}
The distance \( \MMD_D \) metrizes the
weak topology on \(\cP(\Sph)\); that is, for
\((\mu_n)_{n\geq1}\subset\cP(\Sph)\) and \(\mu\in\cP(\Sph)\),
\[
\mu_n\weakstarto\mu
\quad\Longleftrightarrow\quad
\MMD_D(\mu_n,\mu)\longrightarrow0.
\]
\end{proposition}
\begin{proof}
We use Gaussian factorization of the kernel \( K_D \) given in~\eqref{eq:gaussian-factorization}.
For a finite signed Borel measure \(\rho\) on \(\Sph\), define
the finite signed measure
\(\widetilde{\rho}:=D^{1/2}_\#(w\rho)\) on \(\bbR^d\).
Substituting this into~\eqref{eq:rkhs-inner-u} yields
\[
\norm{U_\rho}_{\cH_D}^2
=\inner{U_\rho}{U_\rho}_{\cH_D}
=\iint_{\bbR^d\times\bbR^d}
g(u,v)\,\widetilde{\rho}(\dd u)\widetilde{\rho}(\dd v).
\]
The Gaussian kernel \(g\) is integrally strictly positive definite
\citep[Proposition~5]{sriperumbudur2011universality}, so this
integral is strictly positive whenever
\(\widetilde{\rho}\neq0\).
Since \(D^{1/2}\) is invertible and \(w\) is strictly positive on
\(\Sph\), the map \(\rho\mapsto\widetilde{\rho}\) is injective.
Thus \(\norm{U_\rho}_{\cH_D}^2>0\) for every \(\rho\neq0\).

In particular, the kernel mean embedding is injective on
\(\cP(\Sph)\), so \(K_D\) is \emph{characteristic}.
Since \(K_D\) is continuous and \(\Sph\) is compact,
\citet[Theorem~7]{simongabriel2023metrizing} implies that
\(\MMD_D \) metrizes the weak topology
on \(\cP(\Sph)\).
\end{proof}

We now use this norm to relate the energy gap to the distance from the minimizer.

\begin{proposition}[Energy gap]
\label{prop:energy-gap}
For every \(\nu\in\cP(\Sph)\),
\begin{equation}
\label{eq:gap}
\cE_D(\nu)-\cE_D(\mu_D)
\geq\frac12\MMD_D(\nu,\mu_D)^2.
\end{equation}
\end{proposition}

\begin{proof}
Writing \( U_\nu = U_{\mu_D} + (U_\nu - U_{\mu_D}) \) and using~\eqref{eq:rkhs-norm-sq}, we obtain
\[
\cE_D(\nu)-\cE_D(\mu_D)
= \inner{U_{\mu_D}}{U_\nu - U_{\mu_D}}_{\cH_D} + \frac12\norm{U_\nu-U_{\mu_D}}_{\cH_D}^2.
\]
By linearity of the kernel mean embedding and the reproducing property, the cross-term is
\[
\inner{U_{\mu_D}}{U_\nu - U_{\mu_D}}_{\cH_D}
= \int_{\Sph} \underbrace{\inner{U_{\mu_D}}{K_D(y, \cdot)}_{\cH_D}}_{U_{\mu_D}(y)} \, (\nu - \mu_D)(\dd y).
\]
Therefore,
\begin{equation}
\label{eq:first-variation}
\begin{aligned}
\cE_D(\nu)-\cE_D(\mu_D)
&=\int_{\Sph}U_{\mu_D}\,\dd(\nu-\mu_D)+\frac12\norm{U_\nu-U_{\mu_D}}_{\cH_D}^2 \\
&= \underbrace{\int_{\Sph} \bigl( U_{\mu_D} - 2\cE_D(\mu_D) \bigr) \dd \nu}_{\ge 0} +\frac{1}{2} \MMD_D(\nu,\mu_D)^2,
\end{aligned}
\end{equation}
where we used
\(\int_{\Sph}U_{\mu_D}\,\dd\mu_D=2\cE_D(\mu_D)\)
and \(\nu(\Sph)=1\).
The integral there is nonnegative by the Euler--Lagrange
criterion~\eqref{eq:el}, which proves the inequality.
\end{proof}

\begin{remark}
If \(\supp\mu_D=\Sph\), as will be proved further in
Corollary~\ref{cor:full-support}, then~\eqref{eq:el} gives
\(U_{\mu_D}=2\cE_D(\mu_D)\) on \( \Sph \).
The integral in~\eqref{eq:gap} therefore vanishes, yielding an exact equality
\[
\cE_D(\nu)-\cE_D(\mu_D)
=\frac12\MMD_D(\nu,\mu_D)^2.
\]
\end{remark}

In particular, for every \(\varepsilon\geq0\),
\[
\cE_D(\nu)-\cE_D(\mu_D)\leq\varepsilon
\quad\Longrightarrow\quad
\MMD_D(\nu,\mu_D)\leq\sqrt{2\varepsilon}.
\]
Together with Proposition~\ref{prop:kernel-norm} and weak
continuity of \(\cE_D\), this gives the following characterization
of convergence to the minimizer.

\begin{corollary}[Convergence of approximate minimizers]
\label{cor:approximate-minimizers}
For a sequence \((\nu_i)_{i\geq1}\subset\cP(\Sph)\),
the following are equivalent:
\[
\cE_D(\nu_i)\longrightarrow\cE_D(\mu_D)
\quad\Longleftrightarrow\quad
\MMD_D(\nu_i,\mu_D)\longrightarrow0
\quad\Longleftrightarrow\quad
\nu_i\weakstarto\mu_D.
\]
\end{corollary}

\begin{proof}
The lower bound in~\eqref{eq:gap} gives the first implication,
and Proposition~\ref{prop:kernel-norm} gives the second.
Conversely, weak convergence to \(\mu_D\) implies convergence
of the energies by weak continuity of \(\cE_D\).
\end{proof}

Consequently, a sequence of probability measures has energies
converging to the minimum if and only if it converges weakly to
\(\mu_D\). These observations will allow us to translate energy
estimates into approximation and stability statements.

\input{section_minimizer}

\section{Energy landscape}

A standard approach to the minimization
problem~\eqref{eq:variational-problem} is to equip
\(\cP(\Sph)\) with a suitable metric and consider the
corresponding gradient flow of \(\cE_D\).
Different choices of metric lead to different descent dynamics.
For the Wasserstein distance \(W_2\), the evolution
satisfies
\begin{equation}
\label{eq:wasserstein-flow}
\partial_t\mu_t
=\diver_{\Sph}\bigl(\mu_t\nabla_{\Sph}U_{\mu_t}\bigr);
\end{equation}
see, for instance, \citet{villani2009optimal,santambrogio2017gradient}. A closely related evolution is the normalized repulsive
self-attention equation studied by \citet{geshkovski2025perspective} and \citet{burger2025analysis},
\begin{equation}
\label{eq:normalized-flow}
\partial_t\mu_t
=\diver_{\Sph}\left(
\frac{\mu_t}{U_{\mu_t}}\nabla_{\Sph}U_{\mu_t}
\right)
=\diver_{\Sph}\bigl(\mu_t\nabla_{\Sph}\log U_{\mu_t}\bigr).
\end{equation}
This evolution is a gradient flow of the same energy in a
modified transport geometry with nonlocal mobility
\(U_\mu^{-1}\) (induced by a modified Wasserstein distance). Both equations are understood in the
distributional sense.

The existence of a unique global minimizer does not by itself
ensure convergence of either flow to it.
We therefore study stationary measures
of \(\cE_D\) and show that  every  stationary measure except for $\mu_D$ is singular
with respect to the surface measure. We also remark that these measures can't be local minimizers with respect to the above Wasserstein-type distances. 
We then study the stability
of \(\mu_D\) under perturbations of the interaction matrix.

\subsection{Stationary measures with respect to Wasserstein flows}
\label{sec:stationary-points-singular}

We say that a probability measure
\(\mu\in\cP(\Sph)\) is \emph{stationary} if
\begin{equation}
\label{eq:stationary}
\nabla_{\Sph}U_\mu(x)=0
\qquad\text{for }\mu\text{-a.e.  }x.
\end{equation}
Indeed, testing the respective stationary equations against
\(U_\mu\) gives
\begin{equation*}
0
=\int_{\Sph}\norm{\nabla_{\Sph}U_\mu(x)}^2\,\mu(\dd x),
\end{equation*}
and,
\begin{equation*}
0
=\int_{\Sph}
  \inner{\nabla_{\Sph}U_\mu(x)}
        {\nabla_{\Sph}\log U_\mu(x)}\,\mu(\dd x)
 =\int_{\Sph}
  \frac{\norm{\nabla_{\Sph}U_\mu(x)}^2}{U_\mu(x)}
  \,\mu(\dd x).
\end{equation*}
Since \(U_\mu\) is strictly positive, either identity implies
\eqref{eq:stationary}. Conversely, \eqref{eq:stationary} implies that both
\(\nabla_{\Sph}U_\mu\) and
\(\nabla_{\Sph}\log U_\mu\) vanish
\(\mu\)-almost everywhere, so \(\mu\) satisfies both
stationary equations~\eqref{eq:wasserstein-flow} and~\eqref{eq:normalized-flow}. Hence the constant curve \( \mu_t = \mu\)
solves either evolution. Smoothness of \( K_D \) on the compact sphere
and the bound \( U_{\mu} \geq e^{-\norm{D}_{\op}}\) make both velocity
fields Lipschitz in space and in the measure argument for \( W_2 \).
The corresponding evoluations are therefore unique, so a solution
initialized at \( \mu \) remains there.
This characterization also agrees with stationarity under
transport variations; see
\citet[Lemma~4.3]{burger2025analysis}.

The potential \(U_{\mu_D}\) is constant on \(\Sph\), so
\(\mu_D\) is stationary.
Other stationary
measures exist. For example,
if \(v\in\Sph\) is an eigenvector of \(D\), then
\begin{equation*}
\nabla_{\Sph}U_{\delta_v}(v)
=e^{v^\top Dv}(I-vv^\top)Dv=0,
\end{equation*}
so \(\delta_v\) is stationary. We now show that every stationary measure distinct from
\(\mu_D\) has support of zero surface measure. Thus \( \mu_{D} \) is the only stationary probability measure
having a density with respect to \(\dd \sigma_{d-1}\).

\begin{proposition}[Singularity of non-minimizing stationary measures]
\label{prop:stationary}
Let \(\mu\in\cP(\Sph)\) satisfy \eqref{eq:stationary}.
If \(\mu\neq\mu_D\), then \(\supp\mu\) is contained in the
zero set of a nonzero real-analytic function on \(\Sph\).
In particular,
\begin{equation*}
\sigma_{d-1}(\supp\mu)=0,
\qquad
\mu\perp \dd \sigma_{d-1}.
\end{equation*}
Consequently, \(\mu_D\) is the only stationary probability
measure absolutely continuous with respect to
\(\dd \sigma_{d-1}\).
For \(d=2\), every stationary measure distinct from
\(\mu_D\) has finite support.
\end{proposition}

\begin{proof}
We first verify that \(U_\mu\) extends to a real-analytic
function on \(\bbR^d\). For \(x_0,h\in\bbR^d\), expansion
of the exponential gives
\begin{equation*}
U_\mu(x_0+h)
=\sum_{m=0}^{\infty}\frac1{m!}
 \int_{\Sph}e^{x_0^\top Dy}(h^\top Dy)^m\,\mu(\dd y).
\end{equation*}
The exchange of the sum and the integral is justified by
\begin{equation*}
\frac{e^{x_0^\top Dy}|h^\top Dy|^m}{m!}
\leq
e^{\norm{x_0}\norm{D}_{\op}}
\frac{(\norm{h}\norm{D}_{\op})^m}{m!},
\qquad y\in\Sph.
\end{equation*}
Thus the power series converges absolutely and uniformly
for \(h\) in bounded sets, proving analyticity.

On the sphere,
\begin{equation*}
\nabla_{\Sph}U_\mu(x)
=(I-xx^\top)\nabla U_\mu(x),
\qquad x\in\Sph.
\end{equation*}
It follows that \( g: \Sph \to \bbR \) given by
\begin{equation*}
g(x)
:=\norm{\nabla_{\Sph}U_\mu(x)}^2,
\qquad x\in\Sph,
\end{equation*}
is real-analytic. By \eqref{eq:stationary}, \(g=0\)
\(\mu\)-almost everywhere. Continuity of $g$ then gives
\[
\supp\mu\subset\{x\in\Sph: \; g(x)=0\} =: \cZ(g).
\]
If \(g \equiv 0\), connectedness of \(\Sph\)
implies that \(U_\mu\) is constant on \(\Sph\).
Proposition~\ref{prop:euler-lagrange} then gives
\(\mu=\mu_D\).

Suppose now that \(\mu\neq\mu_D\).
Then \( g \not\equiv 0\).
Since \( \Sph \) is a connected real-analytic manifold,
\citet[Theorem~A.4]{langrenez2024set} gives
\[
\sigma_{d-1}(\cZ(g)) = 0.
\]
Hence \(\sigma_{d-1}(\supp\mu)=0\), which proves singularity.
Since \(\mu_D\) has a smooth, strictly positive
density by Theorem~\ref{thm:mu-d-smooth-density}, it is
the only absolutely continuous stationary probability measure.

For \(d=2\), the zeros of a nonzero real-analytic function
on the circle are isolated by the one-dimensional identity theorem;
see also \citet[proof of Theorem~3.1]{alvarez2026perceptrons}.
Compactness of \( S^1 \) therefore makes the zero set finite,
and the same holds for \(\supp\mu\).
\end{proof}

\begin{remark}
The fact that absolute continuity is a distinguishing property of the minimizer among all stationary points has an important consequence.
If an evolutionary process (e.g., a gradient flow) converges weakly to an absolutely
continuous stationary measure, then its limit must be \( \mu_D \).
In particular, a uniform bound on the relative entropy with respect to
\(\sigma_{d-1}\) ensures absolute continuity of every weak limit,
by lower semicontinuity of relative entropy.
Note, however, that preservation of absolute continuity at each finite time alone
does not suffice to exclude a singular limit.
\end{remark}

\begin{remark}[Absence of local minimizers]

Although the associated Wasserstein-type
gradient flows admit non-minimizing stationary measures, there are no local
minimizers except for \(\mu_D\). Indeed, 
take any metric \(d\) that metrizes weak convergence of measures (e.g., one of the Wasserstein-type distances considered in the previous section) and define 
an open ball
\[
\cV_r(\mu)
=
\{\nu\in\cP(\Sph):d(\mu,\nu)<r\},
\qquad \text{$\mu \in \cP(\Sph)$,  $r>0$}.
\]
Assume \(\nu \in \cP(\Sph)\), $\nu \neq \mu_D$, is a local minimizer of \(\cE_D\) with respect to $d$,
that is, there exists an $r>0$ such that
\[
\cE_D(\nu)\leq\cE_D(\mu)
\qquad\text{for every }\mu\in \cV_r(\nu).
\]
Consider the segment
\[
\nu_t=(1-t)\nu+t\mu_D,
\qquad t\in[0,1].
\]
The map \(t\mapsto\nu_t\) is continuous with respect to weak convergence of measures,
and \(\nu_0=\nu\). Thus, an initial segment is contained in
\(\cV_r(\nu)\): there exists \(t_0\in(0,1]\) such that
\( \nu_{t} \in \cV_r(\nu) \) for every \(t\in[0,t_0)\).
By local minimality and strict convexity~\eqref{eq:convexity}, we have
\[
\cE_D(\nu)
\leq \cE_D(\nu_t)
< (1-t)\cE_D(\nu)+t\cE_D(\mu_D),
\qquad 0<t<t_0.
\]
Subtracting \((1-t)\cE_D(\nu)\) and dividing by \(t\)
gives \(\cE_D(\nu) < \cE_D(\mu_D)\), which contradicts the global minimality of  \(\mu_D\), hence \(\nu=\mu_D\).
\end{remark}

\subsection{Stability with respect to the interaction matrix}
\label{sec:stability}

We now study the dependence of \( \mu_D \) on the interaction
matrix \( D \). The lower bound in~\eqref{eq:gap} allows us to
control the distance between two minimizers by comparing their
energies udner the two interaction kernels.

\begin{proposition}[Stability with respect to the interaction matrix]
\label{prop:parameter-stability}
Let \(E_{1},E_{2}\succ0\) be symmetric, and suppose that
\[
\max\bigl(\norm{E_{1}}_{\op},\norm{E_{2}}_{\op}\bigr)\leq R.
\]
Then
\begin{equation}
\label{eq:parameter-stability}
\MMD_{E_1}(\mu_{E_1},\mu_{E_2})^2
\leq 2e^R\norm{E_1-E_2}_{\op}.
\end{equation}
In particular, if \( D \succ 0 \) and \((D_i)_{i\geq1}\) is a sequence of symmetric
positive definite matrices such that
\(\norm{D_i-D}_{\op}\to0\), then \( \mu_{D_i}\weakstarto\mu_D \).
\end{proposition}

\begin{proof}
For every \(x,y\in\Sph\), the mean value theorem gives
\begin{equation*}
\abs{K_{E_{1}}(x,y)-K_{E_{2}}(x,y)}
\leq e^R\abs{x^\top(E_{1}-E_{2})y}
\leq e^R\norm{E_{1}-E_{2}}_{\op}.
\end{equation*}
Integrating against any \(\rho\in\cP(\Sph)\), we obtain
\begin{equation*}
\abs{\cE_{E_{1}}(\rho)-\cE_{E_{2}}(\rho)}
\leq \frac{1}{2}e^R\norm{E_{1}-E_{2}}_{\op}.
\end{equation*}
By minimality of \(\mu_{E_{2}}\) for \(\cE_{E_{2}}\),
\begin{align*}
\cE_{E_{1}}(\mu_{E_{2}})-\cE_{E_{1}}(\mu_{E_{1}})
&=\bigl(\cE_{E_{1}}(\mu_{E_{2}})-\cE_{E_{2}}(\mu_{E_{2}})\bigr)
+\bigl(\cE_{E_{2}}(\mu_{E_{2}})-\cE_{E_{2}}(\mu_{E_{1}})\bigr) \\
&\quad+\bigl(\cE_{E_{2}}(\mu_{E_{1}})-\cE_{E_{1}}(\mu_{E_{1}})\bigr)
\nonumber\\
&\leq e^R\norm{{E_{1}}-{E_{2}}}_{\op}.
\end{align*}
Applying the lower bound in~\eqref{eq:gap} proves
\eqref{eq:parameter-stability}.

Now suppose that \(\norm{D_i-D}_{\op}\to0\).
The sequence \((D_i)_{i\geq1}\) is bounded in operator norm,
so we may choose \(R\) such that
\(\norm{D}_{\op}\leq R\) and \(\norm{D_i}_{\op}\leq R\)
for every \(i\). Therefore,
\[
\MMD_D(\mu_{D_i},\mu_D)^2
\leq 2e^R\norm{D_i-D}_{\op}
\longrightarrow0.
\]
Proposition~\ref{prop:kernel-norm} gives the claimed
weak convergence.
\end{proof}

\section{Numerical methods}

In this section we use the derived representations of the minimizer
\(\mu_D\) to construct a sampling scheme in arbitrary
dimension and, in dimension two, to approximate the unnormalized
density \(w_D\) by an expansion in Mathieu functions. The code for reproducing the numerical experiments is available at \href{https://github.com/TimRoith/EquilibriumTransformers}{github.com/TimRoith/EquilibriumTransformers}.

\input{section_sampling}
\input{section_mathieu}

\section*{Acknowledgments}
The authors thank the organizers of the 2nd Workshop on Machine Learning in Infinite Dimensions, held at ETH Z\"{u}rich on 8–11 September 2025, where this work originated. Y.K. is grateful to the University of Genova for supporting his research stay in Autumn 2026, when part of this work was undertaken.

\medskip \noindent
\sloppy O.K. acknowledges financial support of the European Commission (Horizon Europe grant ELIAS~101120237).
M.B., S.K. and T.R. acknowledge support from DESY (Hamburg, Germany), a member of the Helmholtz Association HGF. This work was carried in part while T.R. was a member of DESY and was supported in part through the Maxwell computational resources operated at Deutsches Elektronen-Synchrotron DESY, Hamburg, Germany.
M.B. and T.R. acknowledge funding by the German Federal Ministry of Research, Technology and Space (BMFTR) under grant agreement No. 16IS24072B and 16IS24072D (COMFORT).
M.B. and S.K. acknowledge support from the German Research Foundation, project BU 2327/19-2.
T.R. acknowledges the support of COST Action CA24136 - InterCoML, the Munich Center for Machine Learning and the ERC Advanced Grant NEITALG, grant agreement No. 101198055. Funded by the European Union. 
Views and opinions expressed are however those of the authors only and do not necessarily reflect those of the European Union or the European Research Council Executive Agency. 
Neither the European Union nor the granting authority can be held responsible for them.

\begin{center}
\includegraphics[height=1.5cm]{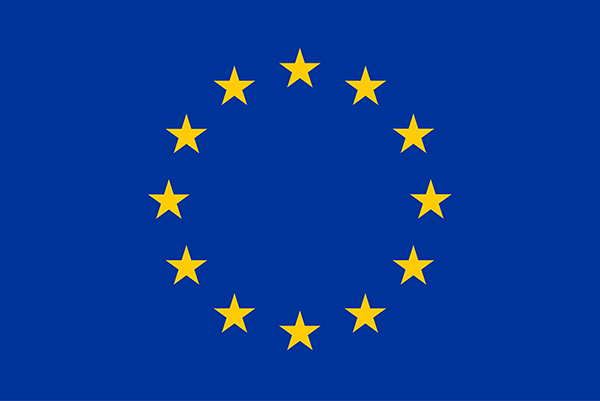}\hspace{1em}%
\includegraphics[height=1.5cm]{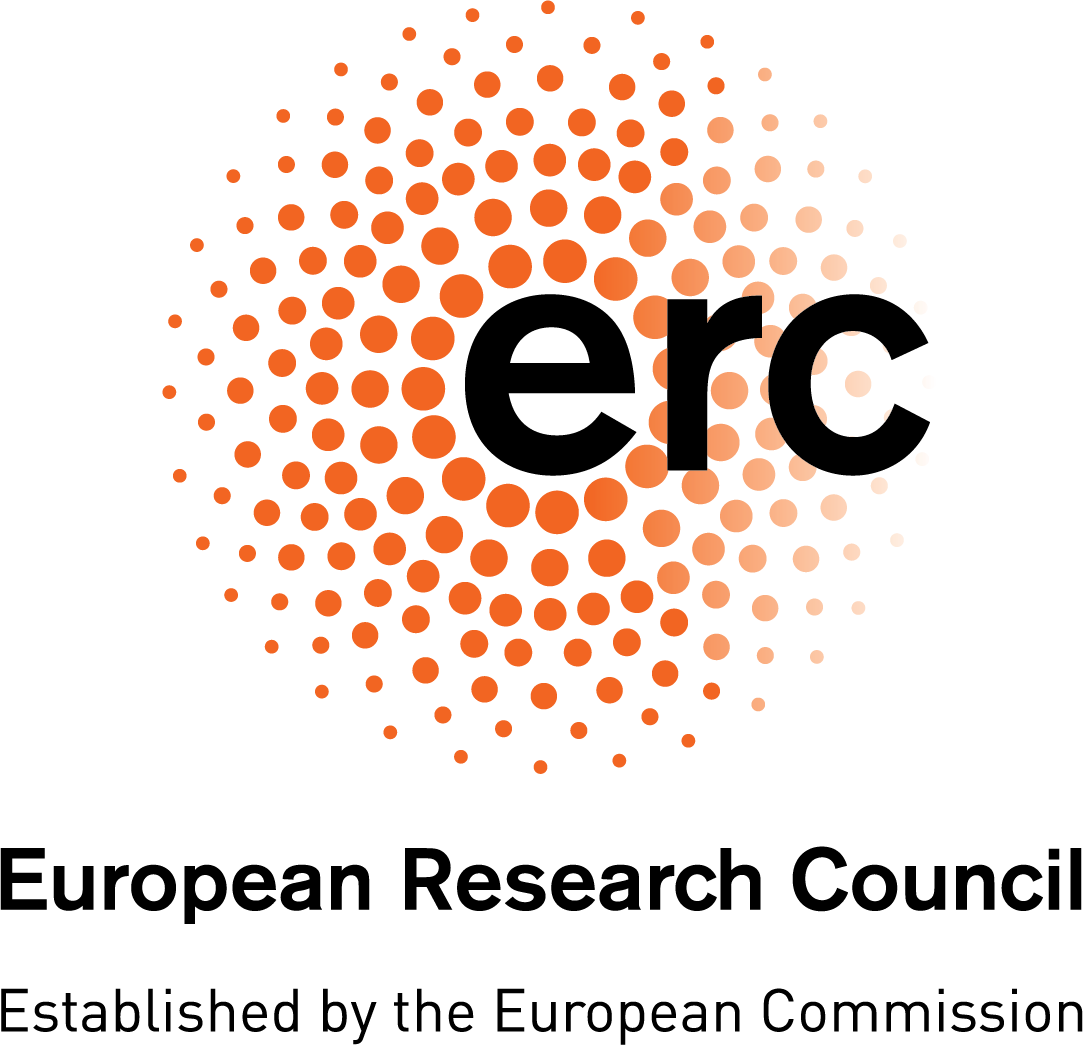}
\end{center}%

\paragraph{Tool and computational resource disclosure.}
Generative AI was used for literature review, assistance with code for the numerical experiments, as well as language and grammar corrections.

\bibliographystyle{plainnat}
\bibliography{references}

\end{document}

%% file: section_minimizer.tex
\section{Properties of the minimizer}
\label{sec:equilibrium}

We are now ready to state the main results. We show that \( \mu_D \) is the distribution of the exit point
of a certain stochastic process and use this representation to establish the existence of
a smooth, strictly positive density on \( \Sph \), resolving conjecture of \citep{burger2025analysis}. This representation also inspires  the numerical sampling
method discussed in Section~\ref{subsec:comp-wos}.

\subsection{A probabilistic representation}
\label{sec:probability}

We start with an informal sketch of the construction; a precise
formulation follows below. Consider the ellipsoid \(D\Ball\),
obtained by applying \(D\) to the unit ball, and start an
isotropic Brownian motion \(B_t = \sqrt{2} W_{t}\) at the origin. Let \(\tau_D\)
denote its \emph{exit time}, that is, the first time the path
reaches the boundary of the ellipsoid, and let \(B_{\tau_D}\) be its exit
location. Independently of the Brownian motion, assign to
each path an exponentially distributed \emph{killing time},
at which the path is deleted. Retain the exit locations
of those paths which reach the boundary before being killed.
Conditional on reaching the boundary before the killing time, the
exit location, mapped back to the sphere by \(D^{-1}\), has law
\(\mu_D\) (Theorem~\ref{thm:cond-exit-measure-minimizer}).
In Section~\ref{sec:density} we use this representation to show
that \(\mu_D\) has a smooth, strictly positive density.

While passing from the sphere to the boundary of an ellipsoid and
back gives a clean geometric intuition,
for the formal argument, we work directly in the unit ball with \( \Sph \) boundary.
Instead of varying the domain, we use an anisotropic Brownian motion, with
different rates of diffusion along the eigendirections of
\(D\).

\paragraph{Tools and notation.}
Let \((\Omega,\cF,\bbP)\) be a probability space supporting a
standard Brownian motion \(W=(W_t)_{t\geq0}\) in \(\bbR^d\),
started at the origin, and write \(\bbE\) for expectation under
\(\bbP\). Here a \(\sigma\)-algebra on \(\Omega\) is a collection
of subsets containing \(\Omega\) and closed under complements
and countable unions. A filtration is an increasing family of
sub-\(\sigma\)-algebras of \(\cF\). We denote by
\(\cF_t:=\sigma(W_s:0\leq s\leq t)\) the natural filtration of
\(W\), where \(\sigma(W_s:0\leq s\leq t)\) is the smallest
\(\sigma\)-algebra making all \(W_s\), \(0\leq s\leq t\),
measurable. A process \((X_t)_{t\geq0}\) is adapted to
\((\cF_t)_{t\geq0}\) if \(X_t\) is \(\cF_t\)-measurable for
every \(t\geq0\). An \emph{event} is a
measurable set \(A\in\cF\). For an introduction to the theoretical tools used in the proof, such as
stochastic processes, martingales, Brownian motion, and stochastic differential
equations (SDE), we refer the reader to
\cite{Baldi2017Stochastic,Williams1991,LeGall2016}.

\begin{theorem}[Conditional exit measure]
\label{thm:cond-exit-measure-minimizer}
Consider a stochastic process \( X \) on \((\Omega,\cF,\bbP)\) given by \( X_t:=\sqrt{2}\,D^{-1}W_t \).
Define
\[
\tau_D:=\inf\{t\geq0:X_t\notin\Ball\}.
\]
Let \(\zeta\sim\Exp(1)\) be independent of the Brownian path, so that \( \bbP(\zeta>t)=e^{-t} \) for \( t \ge 0 \). Then,
\[
\mu_D(A)
=\bbP\bigl(X_{\tau_D}\in A \bigm| \tau_D<\zeta\bigr),
\qquad A\subseteq\Sph \text{ Borel}.
\]
\end{theorem}

\begin{proof}
Let \(\cN_{\bbP}\) denote the collection of all subsets
of \(\bbP\)-null sets in \(\cF\), and set
\[
\cG_t:=\sigma(\cF_t\cup\cN_{\bbP}),
\; t\geq0,
\qquad
\cG_\infty:=\sigma\Bigl(\bigcup_{t\geq0}\cG_t\Bigr).
\]
Then \((\cG_t)_{t\geq0}\) is the completed natural filtration
of \(W\). Since \( D \) is invertible,  \(X\) and \(W\) generate the same natural filtration.
Thus \((\cG_t)_{t\geq0}\) is also the completed natural
filtration of \(X\).
Since \( D \succ 0 \), the first coordinate of
\( X \) is a non-degenerate real Brownian motion, which leaves \( [-1,1] \) almost surely~\citep[Proposition~2.14]{LeGall2016}, therefore \( \bbP(\tau_D<\infty) = 1 \).
\( X_0=0\in\Ball \) and \( X \) has continuous paths, so
\begin{equation}
\label{eq:td-finite-as}
X_{\tau_D}\in\Sph
\; \text{on } \{\tau_D<\infty\}.
\end{equation}
The process \( X \) is adapted to \( (\cG_t)_{t\ge0} \), therefore \( \tau_D \) is a
\( (\cG_t)_{t\ge0} \)-stopping time.

Define event \( E:=\{\tau_D<\zeta\} \).
Conditioning on \( W \) and using the independence of \( \zeta \), we get
\[
M_D:=\bbP(E)=\bbE[e^{-\tau_D}]>0,
\]
since \( \tau_D<\infty \) almost surely.
By~\eqref{eq:td-finite-as},
\[
\bbP\bigl(X_{\tau_D}\in\Sph\mid E\bigr)=1.
\]
Let \( \nu_D \) be given by
\[
\nu_D(A)
:=\bbP\bigl(\{X_{\tau_D}\in A\}\cap E\bigr),
\qquad A\subseteq\Sph\ \text{Borel}.
\]
Since \( X_{\tau_D} \) is measurable, \( \nu_D \) is a Borel measure on
\( \Sph \), with \( \nu_D(\Sph)=M_D \). Hence \( \nu_D/M_D \) is a
Borel probability measure on \( \Sph \), and
\[
\bbP(X_{\tau_D}\in A\mid E)
=\frac{\bbP(\{X_{\tau_D}\in A\}\cap E)}{\bbP(E)}
=\frac{\nu_D(A)}{M_D},
\qquad A\subseteq\Sph\ \text{Borel}.
\]

By Proposition~\ref{prop:unique-minimizer} and the linearity of
\( \rho\mapsto U_\rho \), it is enough to show that \( U_{\nu_D} \)
is constant on \( \Sph \).
Since \( \tau_D \) and \( X_{\tau_D} \) are \( \cG_\infty \)-measurable, conditioning on \( \cG_\infty \) gives
\[
\begin{aligned}
\nu_D(A)
&=\bbE\!\left[
\one_{\{X_{\tau_D}\in A\}}
\one_{\{\tau_D<\zeta\}}
\right]
=\bbE\!\left[
\one_{\{X_{\tau_D}\in A\}}
\bbP(\zeta>\tau_D\mid \cG_\infty)
\right]
=\bbE\!\left[
e^{-\tau_D}\one_{\{X_{\tau_D}\in A\}}
\right],
\end{aligned}
\]
where the last equality follows from the independence of \( \zeta \)
and \( \cG_\infty \), together with \( \bbP(\zeta>t)=e^{-t} \).
Hence, by the standard machine \citep[Section~5.12]{Williams1991},
for every bounded Borel \( \varphi \) on \( \Sph \),
\begin{equation}
\label{eq:nu-d-def}
\int_{\Sph}\varphi\,\dd\nu_D
=\bbE[e^{-\tau_D}\varphi(X_{\tau_D})].
\end{equation}

Fix \( x \in \Sph \). Taking \( \varphi=K_D(x,\cdot) \) gives
\begin{equation}
\label{eq:opt-stopping-uvd}
U_{\nu_D}(x)
= \int_{\Sph} K_D(x, y) \dd \nu_D(y)
=\bbE\!\left[e^{-\tau_D}K_D(x,X_{\tau_D})\right]
=\bbE\!\left[e^{-\tau_D}e^{x^\top D X_{\tau_D}}\right].
\end{equation}
Set \( \beta_t:=x^\top W_t \). Since
\( \norm{x}=1 \), \( \beta \) is a real Brownian motion, and since
\( \beta_t-\beta_s=x^\top(W_t-W_s) \) is independent of \( \cG_s \) for
\( s<t \), it is a \( (\cG_t)_{t \ge 0} \)-Brownian motion. Hence
\begin{equation}\label{eq:zt-exp-martingale}
Z_t
:=\exp\bigl(\sqrt2\,\beta_t-t\bigr)
=e^{-t}e^{x^\top D X_t}
\end{equation}
is a \( (\cG_t)_{t \ge 0} \)-martingale; see, e.g.,
\citet[Theorem~7.5.6]{Durrett2019},
\citet[Section~3.3]{LeGall2016}.
Since \( X_0=0\in\Ball \) and \( X \) has continuous paths,
\( X_{t\wedge\tau_D}\in\overline{\Ball} \) for all \( t\ge0 \), so
\[
0<Z_{t\wedge\tau_D}
=e^{-(t\wedge\tau_D)}e^{x^\top D X_{t\wedge\tau_D}}
\le e^{\norm{D}_{\op}},
\qquad t\ge0.
\]
Thus the stopped martingale \( (Z_{t\wedge\tau_D})_{t\ge0} \) is bounded,
hence uniformly integrable.
Applying the optional stopping theorem~\citep[Theorem~3.22]{LeGall2016}
to this martingale gives
\begin{equation}
\label{eq:optional-stopping}
\bbE\!\left[e^{-\tau_D}e^{x^\top D X_{\tau_D}}\right]
= \bbE Z_{\tau_D}
= \bbE Z_0
= 1.
\end{equation}
Putting~\eqref{eq:optional-stopping} into~\eqref{eq:opt-stopping-uvd}
yields \( U_{\nu_D}(x)=1 \), independent of \( x \). Consequently, by linearity,
\begin{equation}
\label{eq:normalized-potential}
U_{\nu_D/M_D}(x)=\frac{1}{M_D},
\qquad x\in\Sph.
\end{equation}
The claim follows by Proposition~\ref{prop:unique-minimizer}.
\end{proof}

\paragraph{An exit representation of the attention weights.}
Theorem~\ref{thm:cond-exit-measure-minimizer} identifies the
distribution of tokens at equilibrium. We now ask what weight each individual
token assigns to the tokens in its \emph{layer} (or, in the continuous-depth model, at a fixed depth \(t\)). These weights determine the
\emph{attention output} and, through its tangential component, the velocity
of the token in the continuous-depth model.

Let \(\mu\in\cP(\Sph)\) be the distribution of tokens in a given
layer, with each token represented by a point on \(\Sph\).
For a token at \(x\), the \emph{attention weights} define a
probability measure \(Q_{D,\mu}(x,\cdot)\) on \(\Sph\), given by
\begin{equation}
\label{eq:attention-kernel}
Q_{D,\mu}(x,\dd y)
:=\frac{K_D(x,y)}{U_\mu(x)}\,\mu(\dd y).
\end{equation}
Thus \(Q_{D,\mu}\) is a Markov transition kernel within one layer:
for a Borel set \(A\subseteq\Sph\), \(Q_{D,\mu}(x,A)\) is the
total weight that the token at \(x\) assigns to tokens in \(A\).
For an empirical distribution \(\mu=N^{-1}\sum_{j=1}^N\delta_{x_j}\),
these weights are the entries of the corresponding row of the
softmax attention matrix \citep[equation~(1)]{sander2022sinkformers}.

Given a \emph{value} matrix \(V\), the value of a token at \(x\) is \(Vx\).
The \emph{attention output} at \(x\) is the mean of (all) token values weighted by \( Q_{D,\mu}(x,\cdot)\); for example,
for the repulsive choice \(V=-D\) considered here, we write
\begin{equation}
\label{eq:attention-output}
a_{D,\mu}(x)
:=\int_{\Sph}Vy\,Q_{D,\mu}(x,\dd y)
=-\int_{\Sph}Dy\,Q_{D,\mu}(x,\dd y).
\end{equation}
In the continuous-depth model
of \citet[equation~(1.3)]{geshkovski2025perspective},
the position \(x_i(t)\) of the \(i\)-th token evolves according to
\[
\dot x_i(t)
= P_x^\perp \, a_{D,\mu_t}\bigl(x_i(t)\bigr),
\]
where \( P_x^\perp :=I-xx^\top \) denotes the orthogonal projection onto \(T_x\Sph\).
Thus the tangential component of the attention output gives the
token velocity. 

At the equilibrium measure \(\mu_D\), we write
\(Q_D:=Q_{D,\mu_D}\) and \( a_D := a_{D,\mu_D} \). The following corollary
shows that, for each token position \(x\), the attention weights
\(Q_D(x,\cdot)\) correspond to the exit distribution of the diffusion in
Theorem~\ref{thm:cond-exit-measure-minimizer} with a constant drift
determined by \(x\).

\begin{corollary}[Exit distribution of the equilibrium attention weights]
\label{cor:equilibrium-attention-exit}
Let \(X_t=\sqrt{2}\,D^{-1}W_t\) be the process in
Theorem~\ref{thm:cond-exit-measure-minimizer}. For \(x\in\Sph\),
set
\[
X_t^x:=X_t+2tD^{-1}x,
\qquad
\tau_D^x:=\inf\big\{t\geq0:X_t^x\notin\Ball\big\}.
\]
Then \(\bbE[\tau_D^x]<\infty\), and the exit point
\(Y_D^x:=X_{\tau_D^x}^x\) satisfies for \( A\subseteq\Sph\ \text{Borel} \),
\begin{align}
Q_D(x,A)
 &=\bbP(Y_D^x\in A),
 \label{eq:equilibrium-attention-exit}\\
a_D(x)
 &=-2\,\bbE[\tau_D^x]\,x.
 \label{eq:equilibrium-attention-moment}
\end{align}
\end{corollary}

\begin{proof}
Fix \(x\in\Sph\) and set \(T_t:=t\wedge\tau_D^x\). Optional
stopping at the bounded stopping time \(T_t\) gives
\(\bbE[W_{T_t}]=0\) \citep[Corollary~3.23]{LeGall2016}.
Consequently,
\[
\bbE[DX_{T_t}^x]
=\bbE[ \sqrt{2} W_{T_t} + 2 T_t x ]
=2\,\bbE[T_t]\,x.
\]
Taking the scalar product with \(x\) and using
\(X_{T_t}^x\in\overline{\Ball}\) gives
\(2\bbE[T_t]\leq\norm{D}_{\op}\). Monotone convergence gives
\(\bbE[T_t]\to\bbE[\tau_D^x]<\infty\).
Moreover, \(X_{T_t}^x\to Y_D^x\) almost surely and
\(\norm{DX_{T_t}^x}\leq\norm{D}_{\op}\). Hence
\[
\bbE[DY_D^x]
=\lim_{t\to\infty}\bbE[DX_{T_t}^x]
=\lim_{t\to\infty}2\,\bbE[T_t]\,x
=2\,\bbE[\tau_D^x]\,x.
\]

To identify the law of \(Y_D^x\), recall the exponential
martingale \(Z_t=e^{-t}e^{x^\top D X_t}\) defined in~\eqref{eq:zt-exp-martingale}. The stopped martingale \(Z_{t\wedge\tau_D}\) is bounded, and
\(\bbE[Z_{\tau_D}]=1\) by~\eqref{eq:optional-stopping}. We may
therefore define a probability measure \(\bbP_x\) by
\(\dd\bbP_x=Z_{\tau_D}\,\dd\bbP\). Since
\(x^\top D X_t=\sqrt{2}\,x^\top W_t\), Girsanov's theorem
and the Brownian characterization
\citep[Theorems~5.22 and~5.12]{LeGall2016} imply that
\[
\widetilde W_t:=W_t-\sqrt{2}\,(t\wedge\tau_D) \, x
\]
is a standard Brownian motion under \(\bbP_x\). Up to its first
exit from \(\Ball\), the process \(X\) thus satisfies
\[
X_t=\sqrt{2}\,D^{-1}\widetilde W_t+2tD^{-1}x.
\]
Consequently, \(X_{\tau_D}\) under \(\bbP_x\) has the same law
as \(Y_D^x\) under \(\bbP\). Writing \(\bbE_x\) for expectation
under \(\bbP_x\), we obtain, for every bounded Borel function
\(\varphi\) on \(\Sph\),
\begin{align*}
\bbE[\varphi(Y_D^x)]
&=\bbE_x[\varphi(X_{\tau_D})]
=\bbE[Z_{\tau_D}\,\varphi(X_{\tau_D})]
=\bbE\!\left[e^{-\tau_D}K_D(x,X_{\tau_D})
              \varphi(X_{\tau_D})\right] \\
&=\int_{\Sph}K_D(x,y)\varphi(y)\,\nu_D(\dd y)
=\int_{\Sph}\varphi(y)\,Q_D(x,\dd y).
\end{align*}
Here the third equality follows from~\eqref{eq:nu-d-def}.
For the last equality we used \(\nu_D=M_D\mu_D\) and
\(U_{\mu_D}(x)=1/M_D\) from~\eqref{eq:normalized-potential}.
This proves~\eqref{eq:equilibrium-attention-exit}, and the
first-moment calculation gives~\eqref{eq:equilibrium-attention-moment}.
\end{proof}

\begin{remark}
When \(D=\lambda I\) with \(\lambda>0\), the minimizing measure
is \(\sigma_{d-1}\), and \(Q_D(x,\cdot)\) is the von Mises--Fisher
distribution with mean direction \(x\) and concentration
\(\lambda\). In this case Brownian scaling reduces
\eqref{eq:equilibrium-attention-exit} to the first-exit
characterization of \citet{Gatto2013}.
\end{remark}

\begin{remark}
The process \(X\) is an anisotropic Brownian motion, with a covariance matrix
\[
\bbE[X_tX_t^\top]=2tD^{-2}.
\]
In particular, if \(v\) is a unit eigenvector of \(D\) with
eigenvalue \(\lambda\), then \(v^\top X_t\) has variance
\(2t/\lambda^2\). Thus diffusion is faster in directions
corresponding to smaller eigenvalues of \(D\).

This identity gives a comparison between interaction matrices.
Let \(E\succ0\) be symmetric and suppose that
\(D^{-2}-E^{-2} \succeq 0 \). Construct
\[
X_t^{D}:=\sqrt{2}\,D^{-1}W_t,
\qquad
X_t^{E}:=\sqrt{2}\,E^{-1}W_t
\]
using the same Brownian path. Then, for every \(t\geq0\),
\[
\norm{X_t^{D}}^2-\norm{X_t^{E}}^2
=2W_t^\top(D^{-2}-E^{-2})W_t\geq0.
\]
Whenever \(X^{E}\) reaches the sphere, \(X^{D}\) has already
reached it. Their exit times therefore satisfy
\(\tau_D\leq\tau_E\) almost surely, and hence
\[
M_D=\bbE[e^{-\tau_D}]
\geq\bbE[e^{-\tau_E}]=M_E.
\]
By \eqref{eq:normalized-potential} and Proposition~\ref{prop:euler-lagrange},
\[
\min_{\mu\in\cP(\Sph)}\cE_D(\mu)
\leq
\min_{\mu\in\cP(\Sph)}\cE_E(\mu).
\]
In particular, the condition holds when \(D\) and \(E\) admit
a common orthonormal eigenbasis and each eigenvalue of \(D\)
is at most the corresponding eigenvalue of \(E\).
The comparison thus follows by coupling the paths, without
computing either minimizing density.
\end{remark}

\subsection{Smoothness and positivity of the density}
\label{sec:density}

In this section we show that the measure \( \mu_D \) admits a density with respect to the surface measure, and obtain an expression via the Dirichlet
Green function of a certain elliptic operator built from the generator of the diffusion process \( X \). The argument combines the Feynman--Kac formula with Green's representation formula.
Using elliptic regularity and Hopf's boundary point lemma we then prove
smoothness and strict positivity of the density, and hence
full support of the minimizer. For a classical treatment of the elliptic theory used in the proof,
we refer to \citet[Chapters~2, 3, and~6]{gilbarg2001elliptic}.

\begin{theorem}
\label{thm:mu-d-smooth-density}
The measure \( \mu_D \) admits a density with respect to the normalized
surface measure \( \sigma_{d-1} \), given by
\[
\mu_D(\dd y)
=\frac{w_D(y)}{Z_D}\,\dd\sigma_{d-1}(y),
\qquad
Z_D:=\int_{\Sph}w_D\,\dd\sigma_{d-1}.
\]
The function \( w_D \) is smooth and strictly positive on \( \Sph \).
\end{theorem}

\begin{proof}
The \emph{generator} \( \cA_D \) of the process \( X \) is defined on
\( C_c^2(\bbR^d) \) by
\[
\cA_D g(z)
:= \lim_{h\downarrow0} \frac{\bbE_z \left[ g(X_h) \right] - g(z)}{h}
= \lim_{h\downarrow0} \frac{\bbE \left[ g(z+\sqrt2D^{-1}W_h) \right] - g(z)}{h}.
\]
By It\^o's formula~\citep[Theorem~3.3.6]{KaratzasShreve1991},
\[
\cA_D g
=\frac12\Tr\!\left(
(\sqrt2D^{-1})(\sqrt2D^{-1})^\top\nabla^2g
\right)
=\Tr(D^{-2}\nabla^2g)
=\nabla\cdot(D^{-2}\nabla g),
\qquad g \in C_c^2(\bbR^d),
\]
where we used \( D^\top=D \).

For \( f\in C^\infty(\Sph) \), let
\( u\in C^2(\overline{\Ball}) \) be the unique solution of
\begin{equation}\label{eq:dirprob}
\begin{cases}
(\cA_D-\Id)u=0 & \text{in }\Ball,\\
u=f & \text{on }\Sph,
\end{cases}
\end{equation}
whose existence follows from
\citet[Theorem~6.14]{gilbarg2001elliptic}.
Since \( \bbE[\tau_D]<\infty \) by
\citet[Lemma~5.7.4]{KaratzasShreve1991},
applying the Feynman--Kac representation
\citep[Proposition~5.7.2]{KaratzasShreve1991}
with operator \( \cA_D \), killing rate \( k=1 \),
and zero source term \( g=0 \) gives
\begin{equation}
\label{eq:fk}
u(0)
=\bbE\!\left[e^{-\tau_D}f(X_{\tau_D})\right].
\end{equation}

Set \( \cL_D:=\Id-\cA_D \).
The operator \( \cL_D \) is uniformly elliptic, since
\[
\xi^\top D^{-2}\xi
=\norm{D^{-1}\xi}^2
\ge\norm{D}_{\op}^{-2}\norm{\xi}^2,
\]
and it is formally self-adjoint because \( D^{-2} \) is symmetric.
On \( \Sph \), the outward unit normal is \( n(y)=y \).
Set \( n_D(y):=D^{-2}n(y) \), and write
\[
\partial_{n_D}v(y)
:=n_D(y)^\top\nabla v(y)
=n(y)^\top D^{-2}\nabla v(y),
\qquad v\in C^1(\overline{\Ball}\setminus\{0\}),
\]
for the conormal derivative.
Let \( G_D \) denote the negative of the Dirichlet Green function
of \( \cL_D \) in \( \Ball \) with pole at the origin, so that
\[
\begin{cases}
\cL_D G_D = -\delta_0 & \text{in }\Ball,\\
G_D = 0 & \text{on }\Sph.
\end{cases}
\]
Its existence and Green's representation formula are shown in, e.g.,
\citet[Section~2.4]{veron2004elliptic}.
Since \( \cL_Du=0 \), this representation gives
\begin{equation}
\label{eq:green-pre-fk}
u(0)
=\int_{\Sph}
f(y)\,\partial_{n_D}G_D(y)\,\dd S(y).
\end{equation}

Combining~\eqref{eq:nu-d-def} with the Feynman--Kac representation~\eqref{eq:fk}
and~\eqref{eq:green-pre-fk} yields, for every \( f\in C^\infty(\Sph) \),
\[
\int_{\Sph}f\,\dd\nu_D
=\bbE\!\left[e^{-\tau_D}f(X_{\tau_D})\right]
=u(0)
=\int_{\Sph}
f(y)\,\partial_{n_D}G_D(y)\,\dd S(y).
\]
Since \( C^\infty(\Sph) \) is dense in \( C(\Sph) \), this identifies
\( \nu_D=w_D\,\dd S \), where \( w_D:=\partial_{n_D}G_D \).
Elliptic regularity away from the pole gives
\( G_D\in C^\infty(\overline{\Ball}\setminus\{0\}) \)
\citep[Theorem~6.19]{gilbarg2001elliptic}.
In particular, \( \partial_{n_D}G_D \in C^\infty(\Sph) \).
Moreover, the maximum principle gives
\( G_D<0 \) in \( \Ball\setminus\{0\} \), while
\( G_D=0 \) on \( \Sph \).
Hopf's boundary lemma, applied near each boundary point, therefore gives
\[
\partial_nG_D(y)>0,
\qquad y\in\Sph;
\]
see \citet[Lemma~3.4]{gilbarg2001elliptic}.
Since \( G_D = 0 \) on \( \Sph \), i.e., constant, its tangential derivatives
vanish there. Decomposing the ordinary gradient into its
tangential and normal components therefore gives
\[
\nabla G_D(y)
=\nabla_{\Sph}G_D(y)+\partial_nG_D(y)\,y
=\partial_nG_D(y)\,y,
\qquad y \in \Sph.
\]
Consequently,
\[
w_D(y)
=y^\top D^{-2}\nabla G_D(y)
=\bigl(y^\top D^{-2}y\bigr)\,\partial_nG_D(y)>0,
\qquad y\in\Sph.
\]

Write \( a_{d-1}:=\int_{\Sph}\dd S \), so that
\( \dd S=a_{d-1}\,\dd\sigma_{d-1} \).
Since \( M_D=\nu_D(\Sph)=\int_{\Sph}w_D\,\dd S \),
Theorem~\ref{thm:cond-exit-measure-minimizer} gives
\[
\mu_D(A)
=\frac{1}{M_D}\int_A w_D(y)\,\dd S(y)
=\frac{a_{d-1}}{M_D}\int_A w_D(y)\,\dd\sigma_{d-1}(y),
\qquad A\subseteq\Sph\ \text{Borel}.
\]
Thus \( \mu_D \) has the strictly positive smooth density
\( a_{d-1}w_D/M_D \) with respect to \( \sigma_{d-1} \),
and \( \supp\mu_D=\Sph \).
\end{proof}

\begin{remark}
The function \( w_D=\partial_{n_D}G_D \) is the \emph{Poisson kernel}
at the origin for the operator \( \cL_D=\Id-\cA_D \) in \( \Ball \),
with respect to surface measure \( \dd S \).
Thus the density of \( \mu_D \) is obtained by normalizing this kernel.
For a general discussion of Poisson kernels and their relation to
Dirichlet problems, we refer to
\citet[Section~2.4]{veron2004elliptic}.
\end{remark}

\paragraph{Resolution of the conjecture}
The following immediate consequence of the result above resolves the conjecture of
\citet{burger2025analysis}.

\begin{corollary}
\label{cor:full-support}
The density of \(\mu_D\) is bounded away from zero on \(\Sph\):
\[
\frac{w_D(y)}{Z_D}\ge c_D>0,
\qquad y\in\Sph.
\]
In particular, \(\supp\mu_D=\Sph\).
\end{corollary}

\begin{proof}
The continuous, strictly positive function \(w_D/Z_D\) attains
a positive minimum on the compact set \(\Sph\).
\end{proof}

We now formulate the following
conjecture concerning the minimizers of \(w_D\).

\begin{conjecture*}
\label{conj:minimum-density}
The minimum of \(w_D\) on \(\Sph\) is attained at every unit
vector belonging to the eigenspace associated with the largest
eigenvalue of \(D\).
\end{conjecture*}

The probabilistic analysis in Section~\ref{sec:probability} provides intuition for this conjecture. Since the law of the
diffusion is invariant under orthogonal transformations within
each eigenspace of \( D \), the density \( w_D \) takes the same
value at all unit vectors in that eigenspace.
However,
despite the apparent simplicity of the statement,
we have not been able to establish it rigorously.

%% file: section_sampling.tex
\subsection{Sampling from the minimizer}
\label{subsec:comp-wos}

Theorem~\ref{thm:cond-exit-measure-minimizer} gives a construction
which allows to sample from the minimizer without evaluating the density. As explained
in Section~\ref{sec:probability}, one may run
\(B_t=\sqrt2W_t\) in \(\Omega_D=D\Ball\), starting at
the origin, and kill it at an independent exponential
time of rate one. Mapping a surviving exit back by
\(D^{-1}\) gives a sample with law \(\mu_D\).
We approximate these exits by a \emph{Walk on Spheres} (WoS)
\citep{muller1956some,elepov1969solution,sawhney2020monte}.

\paragraph{Survival during one step.}
Suppose that the process starts at the center of a ball
of radius \(r\) contained in \(\Omega_D\). Rotational
invariance makes the exit direction uniform and independent
of the exit time.
For the generator \(\Delta\) and
killing rate one, the probability of reaching the boundary
of this ball before killing is
\begin{equation}
\label{eq:comp-wos-survival}
p_d(r)
:=\frac{r^\nu}{2^\nu\Gamma(\nu+1)I_\nu(r)},
\qquad \nu:=\frac{d}{2}-1,
\qquad p_d(0):=1.
\end{equation}
Here \(I_\nu\) denotes the modified Bessel function of
the first kind of order \( \nu \). By \citet[Equation~(15)]
{rasila2018yukawa}, \(p_d(r)\) is the total mass of the
centered-ball panharmonic measure.
The strong Markov property of
the exponentially killed process allows these steps to be repeated.

\begin{remark}[Numerical evaluation]
For numerical evaluation, it is convenient to use the
exponentially scaled modified Bessel function
\[
I_\nu^{\mathrm{sc}}(r):=e^{-r}I_\nu(r).
\]
Multiplying the numerator and denominator in
\eqref{eq:comp-wos-survival} by \(e^{-r}\) gives, for \(r>0\),
\[
p_d(r)
=\frac{r^\nu e^{-r}}
       {2^\nu\Gamma(\nu+1)I_\nu^{\mathrm{sc}}(r)}.
\]
Evaluating \(I_\nu^{\mathrm{sc}}\) directly avoids overflow
due to the exponential growth of \(I_\nu\) for large
arguments. The value \(p_d(0)=1\) is defined by continuity.
\end{remark}

For \(z\in\Omega_D\), the radius
\begin{equation}
\label{eq:comp-wos-radius}
r(z):=\lambda_{\min}(D)\bigl(1-\norm{D^{-1}z}\bigr)
\end{equation}
gives an admissible ball. Indeed, for every
\(b\in\partial\Omega_D\),
\[
1-\norm{D^{-1}z}
\leq\norm{D^{-1}(b-z)}
\leq\lambda_{\min}^{-1}(D)\norm{b-z}.
\]
Thus \(r(z)\leq\operatorname{dist}(z,\partial\Omega_D)\).
We use a dimensionless tolerance \(0<\varepsilon<1\)
on \(1-\norm{D^{-1}z}\), and project to the sphere when
this quantity is at most \(\varepsilon\).

\begin{algorithm}
\caption{One walk-on-spheres trial for an approximate sample from \(\mu_D\).}
\label{alg:wos}
\begin{algorithmic}[1]
\Require Symmetric \(D\succ0\), tolerance \(0<\varepsilon<1\)
\State \(z\gets0\)
\While{true}
  \State \(y\gets D^{-1}z\), \quad \(\delta\gets1-\norm{y}\)
  \If{\(\delta\leq\varepsilon\)}
    \State \Return \(y/\norm{y}\)
  \EndIf
\State \(r\gets\lambda_{\min}(D)\,\delta\)
  \State Draw \(U\) uniformly on \([0,1]\)
  \If{\(U>p_d(r)\)}
    \State \Return \textsc{Killed}
  \EndIf
  \State Draw \(V\sim\sigma_{d-1}\), independently
  \State \(z\gets z+rV\)
\EndWhile
\end{algorithmic}
\end{algorithm}

\paragraph{Convergence as the tolerance vanishes.}
We next show that the law of an accepted output converges
to \(\mu_D\) as the stopping tolerance tends to zero.

\begin{proposition}[Convergence of the sampling law]
\label{prop:comp-wos-convergence}
Let \(D\succ0\) be symmetric. For \(0<\varepsilon<1\),
denote by \(\mu_{D,\varepsilon}\in\cP(\Sph)\) the law
of the output of Algorithm~\ref{alg:wos}, conditional
on the trial being successful. Then
\[
\mu_{D,\varepsilon}\rightharpoonup\mu_D
\qquad\text{as }\varepsilon\downarrow0.
\]
\end{proposition}

\begin{remark}
Convergence of the classical walk on spheres is recalled
in \citet[Lemma~1.1 and the discussion following it]
{kyprianou2018unbiased}.
Related results for walks with killing and for the error
due to the boundary stopping tolerance are given by
\citet[Section~3.2.1]{bossy2010probabilistic} and
\citet{mascagni2003}, respectively. These results do
not state convergence of the conditional sampling law
for the anisotropic domain (i.e., ellipsoids) and stopping rule used here. We
therefore give a full proof.
\end{remark}

\begin{proof}
Throughout the proof, convergence of random variables
is understood to hold almost surely.

Let \(T_j\) denote the time along the Brownian path
corresponding to the \(j\)-th jump of the walk without
killing. Set \(T_0=0\) and, whenever \(T_j<\tau_D\), define
\begin{align}\label{eq:times}
T_{j+1}:=\inf\{t\geq T_j:
  \norm{D(X_t-X_{T_j})}\geq r(X_{T_j})\}.
\end{align}
If \(T_j=\tau_D\), set \(T_{j+1}=\tau_D\).
By the strong Markov property and the local exit law,
the positions \(X_{T_j}\) follow the jump rule of
Algorithm~\ref{alg:wos}, with killing omitted.
The choice of \(r\) ensures that \(T_{j+1}\leq\tau_D\), so
\(T_j\uparrow T_\infty\leq\tau_D<\infty\) almost surely.
Moreover, continuity of \(X\) gives
\[
X_{T_{j+1}}-X_{T_j}\longrightarrow0.
\]
By the definition of
the successive exit times,
\[
\lambda_{\min}(D)(1-\norm{X_{T_j}})
=\norm{D(X_{T_{j+1}}-X_{T_j})}\longrightarrow0.
\]
Thus \(X_{T_\infty}\in\Sph\), and
\(T_\infty=\tau_D\).

Define \( N_\varepsilon \) the number of jumps needed to reach
the \(\varepsilon\)-neighbourhood of the boundary,
along a given trajectory of the walk without killing, i.e.,
\[
N_\varepsilon:=\min\{j\geq0: 1-\norm{X_{T_j}}\leq\varepsilon\},
\qquad
T_\varepsilon:=T_{N_\varepsilon},
\qquad
Y_\varepsilon:=\frac{X_{T_\varepsilon}}{\norm{X_{T_\varepsilon}}}.
\]
Then \(N_\varepsilon<\infty\), and
\[
\norm{Y_{\varepsilon} - X_{\tau_{D}}}
\leq \varepsilon + \norm{ X_{T_{\varepsilon}} - X_{\tau_{D}}}
\longrightarrow0.
\]
So that
\[
T_\varepsilon\uparrow\tau_D,
\qquad Y_\varepsilon\longrightarrow X_{\tau_D}
\qquad\text{as }\varepsilon\downarrow0.
\]
At each step, the exit location is independent of the
exit time, and \(p_d(r)\) is the probability of surviving
until that exit. The strong Markov property of the killed
process therefore implies that the accepted output has
the law of \(Y_\varepsilon\) conditional on
\(T_\varepsilon<\zeta\).
Its acceptance probability satisfies
\[
\bbP(T_\varepsilon<\zeta)
=\bbE[e^{-T_\varepsilon}]
\geq\bbE[e^{-\tau_D}]>0.
\]
Consequently, for every \(h\in C(\Sph)\), dominated
convergence and Theorem~\ref{thm:cond-exit-measure-minimizer} yield
\[
\int_{\Sph}h\,\dd\mu_{D,\varepsilon}
=\frac{\bbE[e^{-T_\varepsilon}h(Y_\varepsilon)]}
       {\bbE[e^{-T_\varepsilon}]}
\longrightarrow
\frac{\bbE[e^{-\tau_D}h(X_{\tau_D})]}
     {\bbE[e^{-\tau_D}]}
=\int_{\Sph}h\,\dd\mu_D,
\]
recovering the claim.
\end{proof}

Note that the acceptance probability converges to
\(M_D=\bbP(\tau_D<\zeta)\). Combining
\eqref{eq:normalized-potential} with \eqref{eq:el},
we obtain the following relation between the normalization,
the sampling cost, and the minimum:
\[
\bbE[\text{number of exact trials per accepted sample}]
=\frac1{M_D}
=2\cE_D(\mu_D).
\]

\paragraph{Sampling error and the energy gap.}
Independent repetitions,
retaining only successful trials, yield a sequence  of independent random variables
 \((X_i)_{i\geq1}\) with common law
\(\mu_{D,\varepsilon}\). By the strong law of large numbers for empirical measures~\citep{varadarajan1958convergence}, for each fixed
\(\varepsilon\),
\begin{equation}
\label{eq:occupation-wos-measure}
\widehat\mu_{D,\varepsilon,N}
:=\frac1N\sum_{i=1}^N\delta_{X_i}
\rightharpoonup \mu_{D,\varepsilon}
\qquad\text{almost surely as }N\to\infty.
\end{equation}
We now compare the empirical energy of the accepted outputs
with the minimum energy.

\begin{proposition}[Sampling error at fixed tolerance]
\label{prop:comp-sampling-error}
Fix \(0<\varepsilon<1\) and \(N\geq1\). Then,
\begin{equation}
\label{eq:comp-sample-mmd}
\bbE\!\left[
  \MMD_D(\widehat\mu_{D,\varepsilon,N},\mu_D)^2
\right]
\leq
\MMD_D(\mu_{D,\varepsilon},\mu_D)^2
+\frac{e^{\norm{D}_{\op}} - 1}{N},
\end{equation}
where \( \widehat\mu_{D,\varepsilon,N} \) is defined in~\eqref{eq:occupation-wos-measure}.
\end{proposition}

\begin{proof}
Independence gives
\begin{align*}
\bbE[\cE_D(\widehat\mu_{D,\varepsilon,N})]
&=\frac1{2N^2}\sum_{i,j=1}^N\bbE[K_D(X_i,X_j)]\\
&=\frac{N(N-1)}{2N^2}
  \iint_{\Sph\times\Sph}K_D(x,y)\,
    \mu_{D,\varepsilon}(\dd x)\mu_{D,\varepsilon}(\dd y)
+\frac1{2N}\int_{\Sph}K_D(x,x)\,
    \mu_{D,\varepsilon}(\dd x)\\
&=\left(1-\frac1N\right)\cE_D(\mu_{D,\varepsilon})
  +\frac1{2N}\int_{\Sph}K_D(x,x)\,
    \mu_{D,\varepsilon}(\dd x).
\end{align*}
Consequently, Proposition~\ref{prop:energy-gap} yields
\begin{align*}
\bbE\!\left[
  \MMD_D(\widehat\mu_{D,\varepsilon,N},\mu_D)^2
\right]
&=2\left(
  \bbE[\cE_D(\widehat\mu_{D,\varepsilon,N})]-\cE_D(\mu_D)
\right)\\
&=\MMD_D(\mu_{D,\varepsilon},\mu_D)^2
+\frac1N\left(
  \int_{\Sph}K_D(x,x)\,\mu_{D,\varepsilon}(\dd x)
  -2\cE_D(\mu_{D,\varepsilon})
\right)\\
&\leq\MMD_D(\mu_{D,\varepsilon},\mu_D)^2
  + \underbrace{\frac{e^{\norm{D}_{\op}}-1}{N}}_{> 0},
\end{align*}
where the last step uses
\(K_D(x,x)=e^{x^\top Dx}\leq e^{\norm{D}_{\op}}\)
on \(\Sph\) and \(2 \cE_D(\mu_{D,\varepsilon})\geq 1\).
\end{proof}

The first term in \eqref{eq:comp-sample-mmd} accounts
for the stopping tolerance, and the second for the finite
number of accepted samples. Since
\(\mu_{D,\varepsilon} \rightharpoonup \mu_D\) as
\(\varepsilon\downarrow0\), the first term tends to zero.
Taking \(N\to\infty\) and \(\varepsilon\downarrow0\)
therefore gives convergence in mean squared \(\MMD_D\),
and hence convergence of the expected energy to its minimum.
Increasing \(N\) at fixed \(\varepsilon\) does not remove
the stopping-tolerance error. An explicit rate in
\(\varepsilon\) would require a quantitative estimate
of that term.

\begin{figure}
\centering
\includegraphics[width=\linewidth]{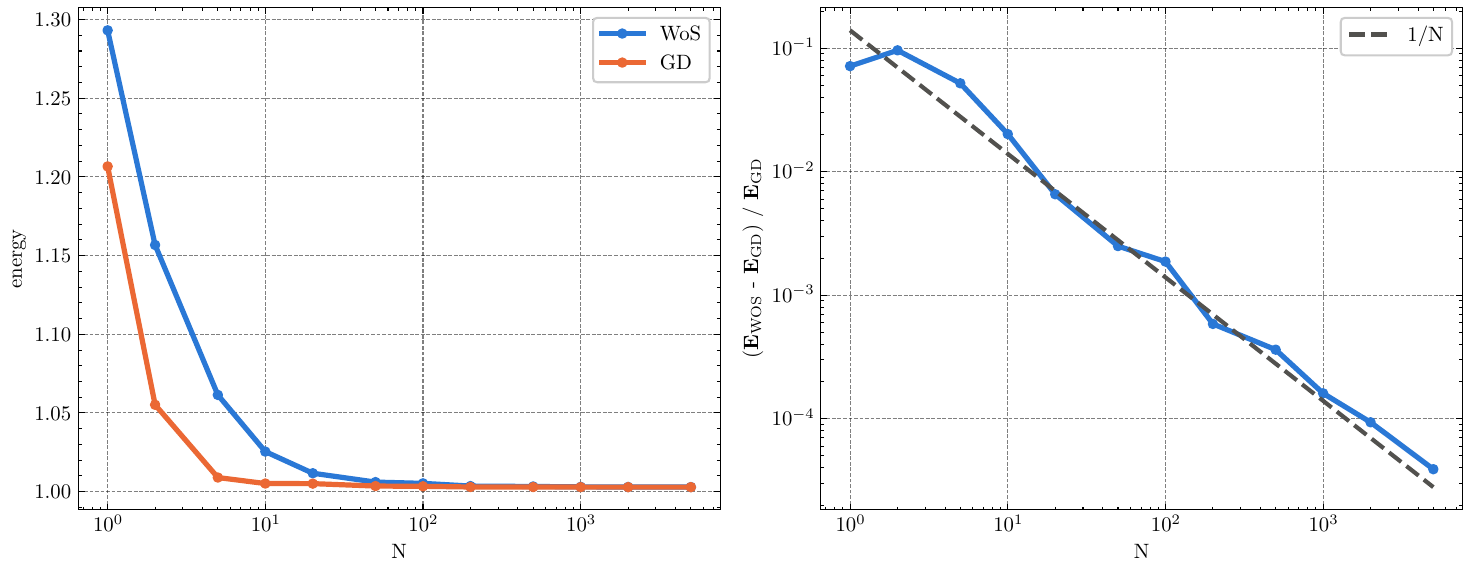}
\caption{Comparison of Algorithm~\ref{alg:wos} with
projected gradient descent \eqref{eq:PGD}, for a fixed
randomly generated matrix \(D\) in dimension \(d=20\).
Left: particle energies at equal ensemble sizes \(N\).
Right: the sampling energy minus the gradient descent
energy, divided by the latter. The dashed line is
proportional to \(N^{-1}\).}
\label{fig:wos}
\end{figure}

For a numerical comparison, we apply the projected gradient
descent scheme from \cite{geshkovski2025perspective}, which optimizes the configuration of an ensemble $\mathbf{x}=(x_1,\ldots, x_N)$. Set $\mathbf{E}_D(\mathbf{x}):=\cE_D(\tfrac{1}{N} \sum_{i=1}^N \delta_{x_i})$. Starting from an initial ensemble $\mathbf{x}^{(0)}$,
we iterate
\begin{align}\label{eq:PGD}
\mathbf{x}^{(k+1)} = 
\Pi
\left(\mathbf{x}^{(k)} - 
\gamma
\nabla_{\mathbf{x}} \mathbf{E}_D(\mathbf{x}^{(k)})
\right),
\end{align}
where $\Pi(z)=z/\|z\|$ for $z\neq0$ is applied to each
particle and $\gamma>0$ is a fixed step size.
We compare the resulting energy with that of the empirical
measure formed from $N$ independent accepted outputs of
Algorithm~\ref{alg:wos}.
Figure~\ref{fig:wos} reports the two energies and their
relative difference for increasing values of \(N\).
Over the range shown, the relative difference decreases
approximately as \(N^{-1}\). %

%% file: section_mathieu.tex
\subsection{Density approximation using Mathieu functions}
\label{sec:Mathieu}

In this section we use the Dirichlet representation to expand
the density $w_D$ of the unnormalized measure $\nu_D$~\eqref{eq:nu-d-def} in an
orthogonal basis. In dimension two, separation of variables
in elliptic coordinates gives an explicit expansion in
Mathieu functions.

Throughout this section, we assume that $D$ is given in diagonal form. Since $\mu_D$ is invariant under coordinate reflections
by~\citet[Corollary~3.10]{burger2025analysis}, its unnormalized
density satisfies $w_D\in L^2_{\text{sym}}(\Sph)$, where \begin{align*}
L^2_{\text{sym}}(\Sph)
:=\bigcap_{r\in\{-1,1\}^d}\left\{f\in L^2(\Sph): f=f\circ \diag(r)\right\},
\end{align*}
is the space of $L^2$ functions on the sphere invariant under coordinate reflections. For every \(r\in\{-1,1\}^d\), the matrix
\(\diag(r)\) is orthogonal and therefore
preserves surface measure on \(\Sph\). Consequently, $f\mapsto f\circ \diag(r)$ defines a linear isometry of \( L^2(\Sph) \). The fixed-point space of each such isometry is a closed
linear subspace of \(L^2(\Sph)\). Hence their intersection,
\(L^2_{\mathrm{sym}}(\Sph)\), is also a closed linear subspace
of \(L^2(\Sph)\). 

For an orthogonal basis $\{f_n\}_{n \in \bbN_0} \subset C^\infty(\Sph)$ of $L^2_{\mathrm{sym}}(\Sph)$
we write the orthogonal expansion of $w_D$ and, for each
integer $K\geq1$, its truncation as
\[
w_{D}
=
\sum_{n \in \bbN_0}
\langle f_n,w_D\rangle\,\frac{f_n}{\|f_n\|^2},
\qquad \text{and} \qquad
w_{D,K}
:=
\sum_{n=0}^K
\langle f_n,w_D\rangle\,\frac{f_n}{\|f_n\|^2}.
\]
Since $\{f_n\}_{n\in\bbN_0}$ is an orthogonal basis of
$L^2_{\text{sym}}(\Sph)$, we have
\[
\lim_{K\to\infty}\|w_{D,K}-w_D\|_{L^2}=0.
\]
In particular, the associated measures converge weakly  \( w_{D,K}\,\dd S \weakstarto \nu_D \) as \( K \to \infty \). For any $n \in \bbN_0$, the Dirichlet problem  \eqref{eq:dirprob} with boundary values $f_n$ admits a unique solution $u_n \in C^2(\overline{\Ball})$. It directly follows from \eqref{eq:green-pre-fk} that \( \langle f_n, w_D \rangle_{L^2} = u_n(0) \). This yields
\begin{align}\label{eq:matseries}
    w_D = \sum_{n \in \bbN_0} u_n(0)\, \frac{f_n}{\|f_n\|^{2}}.
\end{align}

To evaluate the coefficients, we need solutions with a given boundary data. In dimension $d = 2$, pairs of boundary values and solutions are available through Mathieu functions, see \cite{meixner1954mathieusche,whittaker1945course}. These functions arise by separation of variables in the modified Helmholtz equation, 
\begin{align}\label{eq:modifiedhelmholtz}
    \tilde{u} - \Delta \tilde{u} = 0.
\end{align}
More precisely, let us assume that the interaction matrix is given in diagonal form $D = \operatorname{diag}(\lambda_1, \lambda_2)$ with $\lambda_1 > \lambda_2 > 0$. For $c = \sqrt{\lambda_1^2 - \lambda_2^2}$,  we consider the elliptical coordinate transform
\begin{align*}
\Phi: \bbR_{\geq 0} \times (-\pi,\pi] \to \bbR^2,\qquad (\xi, \eta) \mapsto c\, (\cosh(\xi)\cos(\eta), \sinh(\xi)\sin(\eta))^\top.
\end{align*}
We note that for $\xi_0 = \operatorname{arctanh}(\lambda_2/\lambda_1)$ and $\eta \in (-\pi,\pi]$ it holds that
$
    \Phi(\xi_0, \eta) = D (
        \cos(\eta), \sin(\eta)
    )^\top\in D\,\bbS^1.
$ 
We consider solutions of \eqref{eq:modifiedhelmholtz} of the form
\begin{align*}
\tilde{u}_n(\Phi(\xi,\eta))
= \mathrm{Ce}_n(\xi)\,\mathrm{ce}_n(\eta),
\end{align*}
where $\mathrm{ce}_n$ and $\mathrm{Ce}_n$ denote the angular
and modified Mathieu functions, respectively. They satisfy
Mathieu's equation and its modified form
\cite[Sections~28.2, 28.20]{wolf2026}:
\begin{align*}
\text{ce}''_n(\eta) + (a_n - 2 q\cos(2\eta))\cdot \text{ce}_n(\eta) = 0,\qquad
\text{Ce}_n''(\xi) - (a_n - 2q\cosh(2\xi))\cdot \text{Ce}_n(\xi) = 0.
\end{align*}
Here, $q = -c^2/4$ and $a_n \in \bbR$ depends on $q$. Accordingly, the functions $\text{ce}_n, \text{Ce}_n$ depend on $q$, too, however we suppress the dependence on $q$ for brevity.
We refer to \cite[Chapter~2]{mclachlan1951} for a detailed introduction. The functions $\{f_n\}_{n \in \bbN_0}$ defined by
\begin{align*}
    f_n(x) := \text{Ce}_{2n}(\xi_0)\,\text{ce}_ {2n}(\eta),
    \qquad
    x = (\cos(\eta), \sin(\eta))^T,
    \, \eta \in (-\pi,\pi],
\end{align*}
form an orthogonal basis of $L^2_{\text{sym}}(\bbS^{1})$, which can be seen as follows.
As shown in \cite[Chapter~2]{mclachlan1951},
the set $\{\text{ce}_n, \text{se}_{n+1}\}_{n\in\bbN_0}$ is an orthogonal basis of $L^2((-\pi,\pi))$, where $\text{se}_n$ are further solutions of Mathieu's equation, substituting $a_n$ with different values $b_n$.
In fact, $\{\text{ce}_{2n}\}_{n\in\bbN_0}$ consists precisely
of those basis elements satisfying the symmetry conditions
$\text{ce}_{2n}(\eta)=\text{ce}_{2n}(-\eta)$ and
$\text{ce}_{2n}(\pi/2+\eta)=\text{ce}_{2n}(\pi/2-\eta)$,
which correspond to the symmetries defining
$L^2_{\text{sym}}(\bbS^1)$.
In contrast, $\text{se}_{n+1}$ is odd and
$\text{ce}_{2n+1}(\pi/2+\eta)=-\text{ce}_{2n+1}(\pi/2-\eta)$
for every $n\in\bbN_0$.
Hence the corresponding functions on $\bbS^1$ are orthogonal
to $L^2_{\text{sym}}(\bbS^1)$.

We next compute $u_n(0)$ and $\|f_n\|^2$, which determine
the coefficients in \eqref{eq:matseries}.
To this end, let
$\tilde{f}_{2n}:=\tilde{u}_{2n}|_{D\bbS^1}$
denote the restriction of $\tilde{u}_{2n}$ to the ellipse
$D\bbS^1$.
It directly follows that $\tilde{u}_{2n}$ solves the Dirichlet problem
\begin{align*}
        \begin{cases}
        \tilde u - \Delta \tilde u = 0 & \text{in } D\mathbb{B}^2,\\
        \tilde u = \tilde f_{2n} & \text{on } D\bbS^1.
    \end{cases}
\end{align*}
Equivalently, since $f_n = \tilde{f}_{2n} \circ D$, the function $u_n := \tilde{u}_{2n} \circ D$ solves initial Dirichlet problem on $\mathbb{B}^2$ \eqref{eq:dirprob} with boundary value function $f_n$. With the conventional definition of Mathieu functions \cite[Section~2.21]{mclachlan1951},
it holds that $\|\mathrm{ce}_n\|^2 = \pi$. Therefore,
\begin{align*}
    u_{n}(0) = \tilde{u}_{2n}(0) = \text{Ce}_{2n}(0)\, \text{ce}_{2n}(\pi / 2),
\qquad\text{and}\qquad
    \|f_n\|^2 = \text{Ce}_{2n}(\xi_0)^2 \, \|\text{ce}_{2n}\|^2 = \text{Ce}_{2n}(\xi_0)^2\, \pi.
\end{align*}
Note that $\text{Ce}_{2n}(\xi_0)\neq0$: otherwise, the
function $u_n \not= 0$ would solve the Dirichlet problem with
boundary data \( f_n = 0 \), contradicting uniqueness.

Numerical evaluation involves two approximations. We first truncate the series in \eqref{eq:matseries} at order $K\in\bbN$, obtaining 
\begin{align*}
w_{D,K}(x) = \sum_{n=0}^K \frac{\text{Ce}_{2n}(0) \text{ce}_{2n}(\pi/2)}{\text{Ce}_{2n}(\xi_0)\pi} \text{ce}_{2n}(\eta),\qquad \text{for}\quad x=(\cos(\eta), \sin(\eta))^\top. 
\end{align*}%
Since $w_D\in C^\infty(\bbS^1)$, it follows from \cite[Theorem~4.1]{shen2009spectral} that, for every $s\geq 1$
\begin{align*}
\max_{x\in\bbS^1}| w_{D,K}(x) - w_D(x) |
\leq 
C(s,q)\, K^{1/2 - s}\|w_D\|_{H^s}
\end{align*}
where $C(s,q)$ is independent of $K$ and
$\|\cdot\|_{H^s}$ denotes the norm of the Sobolev space
$H^s(\bbS^1)$. In particular, this yields uniform convergence of the approximation as $K\to\infty$.

\begin{figure}
\begin{subfigure}[t]{.5\textwidth}%
\centering%
\includegraphics[width=.9\textwidth]{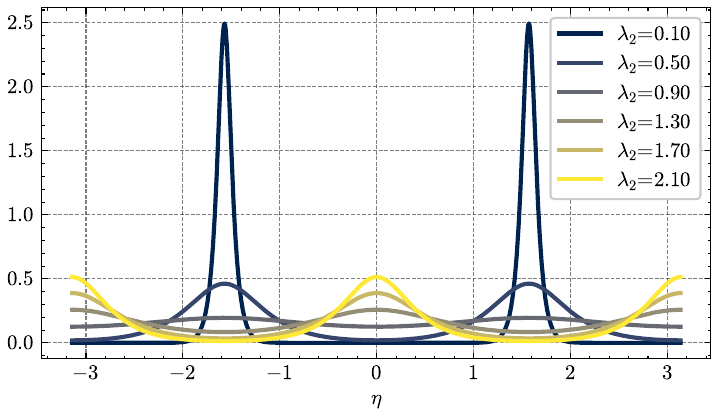}
\caption{Visualization of the approximation $\mathbf{w}_{D,K,N}$.}
\end{subfigure}%
\begin{subfigure}[t]{.5\textwidth}%
\centering%
\includegraphics[width=.9\textwidth]{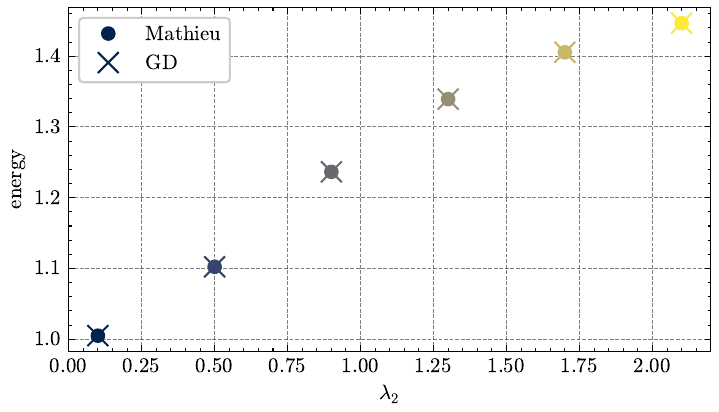}
\caption{Energies of Mathieu approximation compared to the ones obtained with \eqref{eq:PGD} for $N=1000$.
}
\end{subfigure}%
\caption{Mathieu approximation $w_{D,K}$ for matrices $D=\operatorname{diag}(1, \lambda_2)$ with varying $\lambda_2>0$.}\label{fig:math}
\end{figure}

The evaluation of $\text{Ce}_{2n},\text{ce}_{2n}$ is also done approximately by a truncated Fourier series as in \cite[eq. 16--20]{chaos2002mathieu}, where the Fourier coefficients are calculated by the term-recurrences specified therein. 
In Figure~\ref{fig:math} we visualize the approximations obtained this way for different interaction matrices $D$.
We evaluate the energy of an empirical approximation supported
on equally spaced points of $\bbS^1$. For $N\in\bbN$, set
\[
s_\ell=(\cos(\eta_\ell),\sin(\eta_\ell)),
\qquad
\eta_\ell=\frac{2\pi\ell}{N},
\qquad \ell=1,\ldots,N,
\]
and define
\[
\mathbf{w}_{D,K,N}
:= \sum_{\ell=1}^N \left(\frac{w_{D,K}(s_\ell)}{\sum_{j=1}^N w_{D,K}(s_j)} \right) \delta_{s_\ell}.
\]
We then report  $\cE_D(\mathbf{w}_{D,K,N})$.
By Corollary~\ref{cor:full-support}, $w_D$ is bounded away from zero. Uniform convergence therefore
implies the existence of $\mathbf K_D \in \bbN$ such that
$w_{D,K}>0$ on $\bbS^1$ for every $K \geq \mathbf K_D$.
In the implementation, we start with $K=60$ and double $K$
until all evaluated weights are positive. This condition
ensures that the normalized discrete measure is a probability
measure; it does not provide an estimate of the truncation
error.

